\documentclass[11pt,a4paper]{article}
\usepackage[dutch,english]{babel} 
\usepackage{fullpage}
\usepackage[round]{natbib} 
\usepackage{graphicx} 
\usepackage{graphics} 
\usepackage{amsmath} 
\usepackage{amsfonts}
\usepackage{amssymb, amsthm}
\usepackage{bbm}

\usepackage{mathtools}

\usepackage{lscape}

\usepackage{array}
\usepackage{siunitx}
\usepackage{floatrow}
\usepackage{subcaption}
\usepackage{soul}
\usepackage{mathtools}

\usepackage{lmodern}
\usepackage{rotating}
\usepackage[ruled, vlined, linesnumbered]{algorithm2e}

\usepackage{tikz}
\usepackage{tabularx}
\usepackage{multirow}
\usetikzlibrary{automata,positioning}
\usepackage{dcolumn}
\usepackage{subcaption}
\usepackage{enumitem}
\usepackage{abstract}
\usepackage[left=2.2cm,right=2.2cm,top=2.2cm,bottom=2.2cm,twoside]{geometry}
\usepackage{geometry}
\usepackage{setspace}
\usepackage{ragged2e}
\usepackage{float}
\usepackage{colortbl}
\usepackage[flushleft]{threeparttable}
\usepackage{booktabs}
\usepackage{caption}
\usepackage{booktabs,threeparttable}
\usepackage{placeins}
\usepackage{caption}
\usepackage{breqn}
\usepackage{verbatim}
\usepackage{eurosym}
\usepackage{changepage}
\usepackage{listings}
\usepackage{xcolor}
\usepackage{fancyvrb}

\usepackage{color} 
\definecolor{mygreen}{RGB}{28,172,0} 
\definecolor{mylilas}{RGB}{170,55,241}

\newcommand{\tup}[1]{\textup{#1}}

\begin{document}
\LARGE \noindent A Continuum Approximation Approach to the Hub Location Problem in a Crowd-Shipping System

\normalsize

\begin{flushleft}
\noindent Patrick Stokkink$^1$, Nikolas Geroliminis$^1$\\
\scriptsize 
$^1$ \'{E}cole Polytechnique F\'{e}d\'{e}rale de Lausanne (EPFL), Urban Transport Systems Laboratory (LUTS), Switzerland\\
 Contact: patrick.stokkink@epfl.ch (Patrick Stokkink),  nikolas.geroliminis@epfl.ch (Nikolas Geroliminis)
\end{flushleft}

\noindent\rule[0.5ex]{\linewidth}{1pt}

\normalsize \noindent \textbf{Abstract}\\ \small
Last-mile delivery in the logistics chain contributes to emissions and increased congestion. Crowd-shipping is a sustainable and low-cost alternative to traditional delivery, but relies heavily on the availability of occasional couriers. In this work, we propose a hub-based crowd-shipping system that aims to attract sufficient potential crowd-shippers to serve a large portion of the demand for small parcels. While small-scale versions of this problem have been recently addressed, a scaling to larger instances significantly complexifies the problem. A heuristic approach based on continuum approximation is designed to evaluate the quality of a potential set of hub locations. By combining an efficient and accurate approximation method with a large neighborhood search heuristic, we are able to efficiently find a good set of hub locations, even for large scale networks. Furthermore, on top of determining good hub locations, our methods allow to identify the expected number of delivered parcels in every region, which can be used to design a smart dynamic assignment strategy. \\
A case study on the Washington DC network shows that hubs are built at locations that are both geographically central, but most importantly are popular origins for crowd-shippers. The optimal number of hubs is mainly dependent on the marginal number of parcels that can be served by crowd-shippers from a specific hub, relative to the costs involved in opening that hub. The performance of our algorithm is close to that of a simulation-optimization algorithm, yet being up to 25 times faster. Thereby, the results show a dynamic assignment policy based on continuum approximation estimates outperforms existing assignment strategies. 

\\
\noindent\rule[0.5ex]{\linewidth}{1pt}

\section{Introduction} \label{sec:introduction}
\begin{figure}[H]
    \centering
    \includegraphics[width=0.7\textwidth]{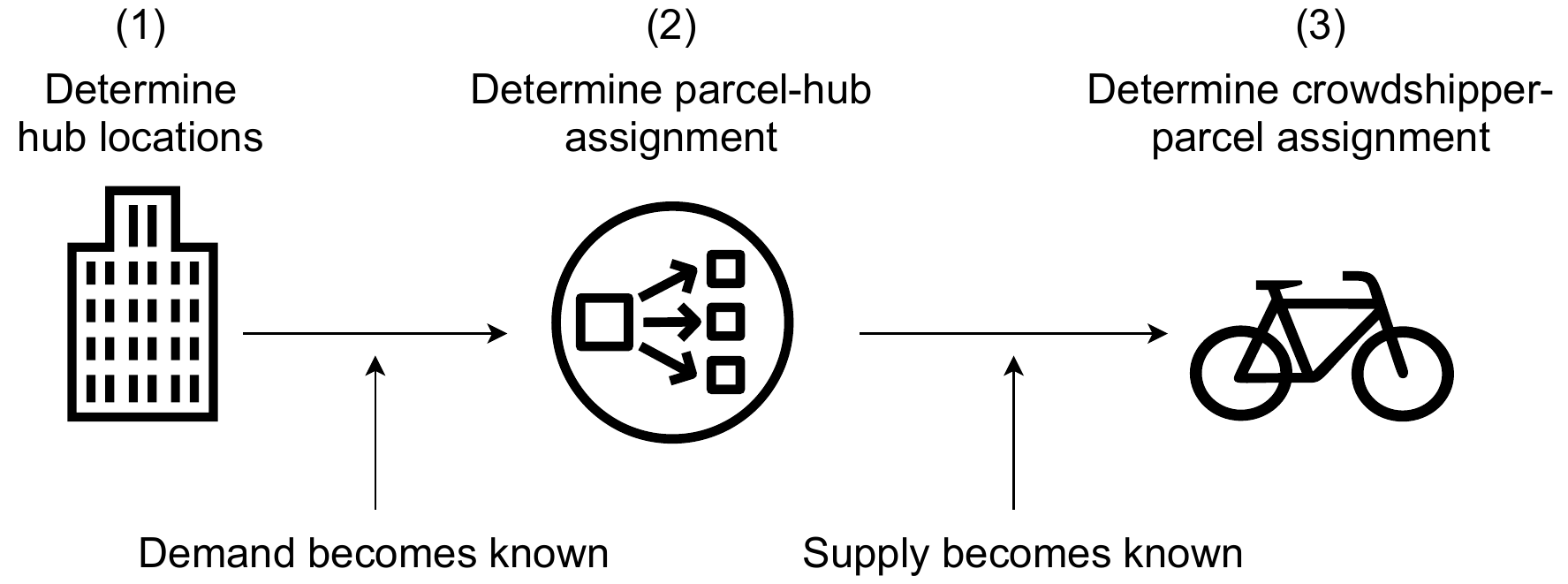}
    \caption{Schematic representation of decision process}
    \label{fig:problem_scheme}
\end{figure}

The growing demand for e-commerce has led to a substantial increase in the challenges faced in traditional delivery. Nowadays, the ``sharing economy" allows to rapidly connect supply and demand, which can be used to overcome these challenges. Such a system where last-mile delivery is outsourced to a large number of individuals is referred to as \textit{crowd-shipping}. In a crowd-shipping system, individual couriers perform deliveries on their pre-existing route, possibly with a small detour, and thereby contribute to the last-mile delivery of small parcels. 
\\\\
Crowd-shipping has numerous advantages for customers, retailers and the society as a whole. Customers, who can be either on the supply or demand side, are offered a fast, flexible and often cheap alternative of delivery. Retailers also benefit in terms of cost, as occasional couriers (hereafter also referred to as crowd-shippers) are generally cheaper than standard drivers. Thereby, these occasional couriers are more flexible compared to standard drivers routes which are planned in advance, and can therefore more easily incorporate last-minute changes. Society benefits mostly from the reduced environmental impacts as well as reduced traffic congestion. Traditional delivery is a large source of congestion because of the presence of large delivery vans in urban networks and the road-blocks they cause. According to \cite{carbonfootprint2020}, the vast majority of the emissions in the whole logistics chain is generated by last-mile delivery. On top of that, polluting delivery vehicles may not be allowed to enter low-emission or zero-emission zones, making it absolutely necessary to look for green alternatives such as crowd-shipping.
\\\\
One of the most important drivers of a crowd-shipping system is the availability of supply and the potential to match supply to demand without large detours of crowd-shippers. The potential pool of crowd-shippers is considered to have a planned personal trip and an associated individual trajectory which is not necessarily near to a parcel trajectory. In this work we focus on few-to-many delivery problems of small parcels that are transportable by foot or by bike. Clearly, if the origins of demand are poorly accessible by potential crowd-shippers, supply is lost and demand remains unserved. For this reason, we focus on the problem of determining optimal hub locations that function as origins for parcels (demand). Good hub locations are built at central locations in the network such that they serve a large portion of demand by attracting sufficient crowd-shippers (supply). 
\\\\
In a hub-based crowd-shipping system, various decisions have to be made to construct a profitable system. These decisions can be divided into strategic, tactical and operational decisions. A schematic representation of the decision process is illustrated in Figure \ref{fig:problem_scheme}. In the first stage, the locations of the hubs (i.e. locations where parcels are stored for crowd-shippers to pick them up) have to be determined without knowledge of the exact demand and supply realizations. This is a strategic decision which has to be made before the system is operational. Thereafter, demand becomes known either partially or completely and in the second stage the assignment of parcels to these hubs are determined. This is a tactical decision which is based on observed demand and expectations of supply, fed by historical data. Then, the supply becomes known, usually in a dynamic fashion, and the parcels are assigned  to crowd-shippers in the third and final stage. These operational decisions are made on a daily basis based on known demand and partially known supply. 

\FloatBarrier
\FloatBarrier

\noindent In this paper, we develop a framework that allows to determine the best hub locations for a crowd-shipping system in a large urban area. This problem is especially difficult because of the dependency of lower-level decisions and costs on upper-level decisions. To track these interactions, we solve the lower-level assignment problem of parcels to potential crowd-shippers through a Continuum Approximation (CA) approach, allowing us to determine the lower-level costs efficiently in short time. These estimates are based on physical properties of the matching procedure, as well as expectations of supply and demand, fed by historical data. We develop a large neighborhood search heuristic that exploits the CA estimates to efficiently search a good set of hub locations that minimizes the operational costs. In addition to this, these estimates are used to design a smart dynamic assignment strategy of parcels to crowd-shippers that outperforms existing strategies. A comparison of our approach to a common simulation-optimization approach shows very similar performance in objective value, whereas the computation time of our approach is substantially lower. 
\\\\
The remainder of this paper is organized as follows. We provide a review of the relevant literature in Section \ref{sec:literature}, where we highlight the unique features of this new problem that require to introduce new formulations and solution approaches compared to similar existing problems in the literature. The methods we propose to efficiently and simultaneously solve this multi-level problem, including the continuum approximation of the assignment problem, the hub-location algorithm and a smart dynamic assignment strategy, are discussed in Section \ref{sec:methodology}. The results are discussed in Section \ref{sec:results}, where we evaluate the proposed methodology with comparisons of benchmark methods and a much slower simulation-based approach, using a discrete event simulator based on a part of the city of Washington DC. The paper is concluded in Section \ref{sec:conclusion}.

\section{Literature Review} \label{sec:literature}

The last-mile delivery of products to people is a well-studied topic in the optimization literature. Traditionally, goods are delivered by using delivery vans. In this case, the problem can be formulated as a Pickup-and-Delivery Problem (PDP) (\cite{savelsbergh1995general}) or a Vehicle Routing Problem (VRP) (\cite{toth2002vehicle}). These problems have been extended to include various problem-specific aspects such as time-windows (e.g. \cite{dumas1991pickup}, \cite{ropke2009branch}) or uncertainty (e.g. \cite{fabri2006dynamic}).
\\\\
Due to the increase of online shopping, more traditional vans are needed to serve all demand. According to \cite{iwan2016analysis}, delivery vans for last-mile delivery are one of the main causes of congestion in urban areas. As a consequence, many companies are looking for more sustainable options to replace the aforementioned traditional delivery methods. \cite{iwan2016analysis} analyze the use of parcel lockers where customers can pickup and send small parcels. Their results from a pilot survey in Poland indicate that the use of these lockers can potentially reduce the environmental impact of last-mile delivery. Another alternative is drone delivery \citep{murray2015flying, agatz2018optimization, karak2019hybrid}, for which it has been shown that combining truck delivery with a drone can significantly reduce transportation cost. Thereby, drones cause less congestion compared to delivery vans.  \cite{akeb2018building} propose a model that relies on the interaction of a network of neighbors to enhance parcel delivery in urban areas.
\\\\
In a crowd-shipping system, the last-mile delivery of small parcels is (partially) outsourced to individual commuters that can deliver the parcel on their pre-existing route. Various empirical studies have investigated the potential and determinants of crowd-shipping (e.g. \cite{marcucci2017connected}, \cite{ermagun2018bid}, \cite{le2018crowd}, \cite{le2019influencing}, \cite{punel2019push}). These studies have shown the potential demand for crowd-shipping and the concerns of potential users. Thereby, they highlight the importance of the availability of supply. Potential crowd-shippers have a pre-existing itinerary (origin, destination and approximate departure and arrival times) and trip-purpose (for example, a work commute or a leisure trip). Therefore, only parcels can be assigned for which the crowd-shipper does not experience too much inconvenience and the inconvenience has to be compensated for through some (monetary) incentive.
\\\\
Recently, substantial research has been done on the operational problems that arise in crowd-shipping systems. For a review on the recent academic research as well as recent practice, the reader is referred to \cite{le2019supply}. \cite{pourrahmani2021crowdshipping} give an overview of the operational challenges and research opportunities that exist in this field. One of these operational challenges is the matching of parcels to crowd-shippers which has been studied by, among others, \cite{li2014share} and \cite{soto2017matching}. We also note the similarities with matching problems in ride-sharing (\cite{masoud2017real}), ride-hailing (\cite{yan2020dynamic} and carpooling (\cite{wang2018stable}). Another important operational problem is pickup and delivery routing. Clearly, these problems are intertwined and therefore these two problems are often tackled jointly. \cite{archetti2016vehicle} model the static routing problem as a Vehicle Routing Problem with Occasional Drivers (VRPOD). It is assumed that an occasional driver is willing to make a delivery if the extra distance travelled to make the delivery is less than a pre-specified portion of the total distance travelled. \cite{li2014share} consider a crowd-shipping scenario where people and parcels share a taxi. They therefore model the problem as a Share-a-Ride Problem (SARP), which is an extension of the Dial-a-Ride Problem (DARP). \cite{dahle2017vehicle} consider a two-stage stochastic program to model the VRP with dynamic occasional drivers. The first-stage decision models the route of the traditional delivery truck. After the occasional drivers make themselves known in the second stage, they are assigned to parcels and the truck route can be changed. \cite{arslan2019crowdsourced} model the problem as a dynamic PDP. Their heuristic assigns crowd-shipping tasks to occasional (ad-hoc) drivers dynamically. \cite{cohn2007integration} integrate matching and routing decisions for carriers of small packages. \cite{yildiz2019service} introduce service and capacity planning problem. With their model, they aim to answer questions that arise in a crowd-shipping system, both on strategic and operational levels.
\\\\
As the availability of supply is a key determinant of the performance of a crowd-shipping system, parcels may be stored at intermediate hub locations (or transshipment points) such that they are easily reachable by potential suppliers. \cite{wang2016towards} consider ``pop-stations" distributed around the city where crowd-shippers can perform pickups. For a fixed set of transshipment points, they optimize the utilization of crowd-shippers for last-mile delivery. \cite{raviv2018crowd} and \cite{macrina2020crowd} consider a crowd-shipping system where crowd-shippers can pickup parcels either from the depot or from transshipment points. Their results show the economic benefits of such transshipment nodes. Similarly, \cite{yildiz2021express} also considers transshipment points but uses a dynamic programming algorithm to solve their problem. Contrary to the fixed transshipment points in the previous works, \cite{mousavi2021stochastic} consider mobile depots. They do not consider the routing of vehicles, but they determine the optimal location of these mobile depots under uncertainty in supply.
\\\\
Contrary to the majority of the literature that have studied operational problems, in this paper we consider the strategic planning problem of network design. Specifically, we focus on finding the optimal hub locations with uncertainty in supply and demand. The problem we study is quite closely related to the hub-location-and-routing problem \citep{aykin1995hub,wasner2004integrated,rodriguez2014branch}. In such problems, the optimal hub locations, corresponding service areas and routes are jointly determined. The most fundamental difference between this problem and our problem is the existence of crowd-shippers. The use of crowd-shippers may cause large inter-dependencies between the performance of various hubs, especially when they are build close to each other. Similar to the hub-location-and-routing problem, in this problem we jointly determine the hub locations, approximate service areas and the assignment of parcels to crowd-shippers. As the tactical and operational decisions in stage 2 and 3 are strongly interconnected with the hub locations, this does not allow us to decompose the problem. 
\\\\
Exact formulations for strategic location problems in crowd-shipping systems have been proposed by \cite{mousavi2021stochastic} and \cite{nietovalue}, who used exact approaches to solve small-scale instances. Given that crowd-shipping systems mostly exist in large urban areas and their surrounding suburbs, the considered networks are of such as size that solving this difficult problem exactly is infeasible. For this reason, we use a CA approach to approximate with some level of accuracy the cost of the lower-level routing and matching decisions and consecutively use an efficient heuristic to determine the hub locations. CA approaches have been widely used for the design of large-scale networks. For example, for the design of vehicle routing problems \citep{daganzo2005logistics, ouyang2007design}, integrated package distribution systems \citep{smilowitz2007continuum} and pickup-and-delivery problems \citep{lei2018continuous}. Thereby, CA has been used for various variants of the Facility Location Problem (FLP), such as the reliable facility location problem \citep{li2010continuum, cui2010reliable} and the competitive facility location problem \citep{wang2013continuum}. Although we note the similarities between the FLP and the hub location problem discussed in this paper, we note that the main difference is that the facility location problem is based on a distance metric and nearest-assignment strategies, whereas our problem introduces an additional layer of complexity due to the assignment problem of parcels to crowd-shippers. Due to this complexity, crowd-shippers cannot be assigned to parcels in a nearest-neighbor approach as this will lead to large suboptimalities.  
\\\\
Building on the approximated lower level costs, we use a large neighborhood search heuristic to solve the hub location problem that minimizes the total cost of the crowd-shipping system. We use performance metrics to efficiently search the neighborhoods of solutions. Due to the fast CA approximation, we are able to find the optimal hub locations in reasonable time. In addition to this, the CA estimates are used as an input to a smart dynamic assignment strategy, that outperforms existing dynamic strategies by leveraging expectations of future crowd-shippers. A discrete event simulator is used to evaluate the performance of the continuum approximation, the hub-location algorithm and the dynamic assignment strategy.

\section{Methodology} \label{sec:methodology}
In this section we describe the methodological contributions of this paper. In Section \ref{subsec:discrete} we give a more detailed description of the problem and formulate the problem as an integer stochastic programming problem. In Section \ref{subsec:ca} we describe the continuum approximation approach used to approximate the cost in the third-stage assignment problem. Section \ref{subsec:stage_1} explains the methods used to identify the hub locations based on the CA results. Section \ref{subsec:stage_2} describes the parcel-hub assignment procedure and Section \ref{subsec:stage_3} describes the (dynamic) assignment procedure of parcels to crowd-shippers, which was previously approximated. A notational glossary of the sets, parameters and variables is provided in Table \ref{tab:notation}.

\begin{table}[H]
\caption{Notational Glossary}
\label{tab:notation}
\footnotesize 
\begin{tabular}{ll}
\toprule 
$a_{rh}$ & Number of parcels with final destination $r$ stored at hub $h$\\
$B$ & Set of potential crowd-shippers in a batch\\
$C_{\tup {total}}(H)$ & Total expected cost of opening the set of hubs $H$\\ 
$D$ & Set of demand requests\\
$d_r$ & Actual demand in region $r$  \\ 
$\hat{d}_r$ & Expected demand in region $r$  \\
$\hat{d}_r'$ & Remaining expected demand in region $r$ in iterative CA procedure \\
$e_{ijhr}$ & Binary parameter indicating if crowd-shippers travelling from $i$ to $j$ can feasibly pick up a parcel at hub $h$\\ & and deliver it to $r$\\
$\Tilde{e}_{ijr}$ & Binary parameter indicating if crowd-shippers travelling from $i$ to $j$ can feasibly pick up a parcel at\\ & \textit{at least one} hub and deliver it to $r$\\
$f$ & Fixed cost of opening a hub\\
$H$ & Set of opened hubs\\
$H^p$ & Set of potential hubs (or hub locations) \\
$l_r$ & Leftover demand in region $r$ in iterative CA procedure\\
$n$ & Number of possible hub locations\\
$o_h$ & Decision variable indicating whether hub $h$ is opened\\
$p_{\tup {cs}}$ & Reward given to crowd-shippers to perform a pickup and delivery\\
$p_{\tup {reg}}$ & Cost of regular delivery of one parcel\\
$Q$ & Maximum number of hubs that can be built\\
$R$ & Set of regions\\
$S$ & Set of potential suppliers \\
$s_{ij}$ & Total demand that can potentially be served by crowd-shippers travelling between $i$ and $j$ \\
$s_{h_1h_2}$ & Similarity of hubs $h_1$ and $h_2$\\
$u_{ds}$ ($u_{dsh}$) & Parameter indicating if crowd-shipper $s$ can feasibly pick up and deliver demand request $d$\\ 
& (when it is assigned to hub $h$)\\
$\Bar{u}_d$ & Number of potential crowd-shippers that can be feasibly matched to a specific demand request $d$ \\
$v_h$ & Quality metric for opening single hub $h$ \\
$w_{dh}$ & Decision variable indicating whether demand request $d$ is assigned to hub $h$\\
$x_{ds}$ & Decision variable indicating whether demand request $d$ is matched to supplier $s$\\
$y_r$ & Intermediate estimate for expected number of parcels deliver bt crowd-shippers to region $r$\\
$z_r$ & Expected number of parcels delivered by crowd-shippers to region $r$  \\
$z_r(H)$ & Total expected number of parcels delivered by crowd-shippers to region $r$ for a set of opened hubs $H$ \\
$z_r^h$ & Expected number of parcels delivered by crowd-shippers to region $r$ when only hub $h$ is opened \\
$\alpha$, $\beta$, $\gamma$ & Tuning parameters\\
$\eta$ & Number of initial solutions in hub-location algorithm\\
$\kappa$ & Iteration limit of hub-location algorithm\\
$\lambda_{ij}$ & Actual number of potential crowd-shippers travelling between regions $i$ and $j$\\
$\hat{\lambda}_{ij}$ & Expected number of potential crowd-shippers travelling between regions $i$ and $j$\\
$\hat{\lambda}_{ij}'$ & Remaining expected crowd-shippers travelling between regions $i$ and $j$ in iterative CA procedure\\
$\xi_D$ ($\xi_S$) & Random variable of which the realizations are vectors of demand requests (potential suppliers)\\
$\tau$ & Maximum detour crowd-shippers are willing to make to pick up and deliver a parcel\\
$\psi_k(\cdot)$ & Function describing the costs of the $k^{th}$ stage for a given set of inputs\\
$\Omega$ & Current solution (i.e. set of hubs) in the hub-location algorithm\\

\bottomrule
\end{tabular}
\end{table}

\subsection{Discrete Formulation}\label{subsec:discrete}
The problem as described in Figure \ref{fig:problem_scheme} can be formulated as a mixed integer stochastic programming problem. In this way, we can incorporate the two types of uncertainty by dividing the problem into three levels. The formulation is similar to that of \cite{mousavi2021stochastic}. The main difference is that we consider an additional layer of uncertainty (uncertainty in demand) which makes our problem a three-stage stochastic programming problem, compared to the two-stage stochastic programming problem proposed by \cite{mousavi2021stochastic}. 
\\\\
We consider a set of demand requests $d \in D$, which is a realization drawn from random variable $\xi_D$. Similarly, we consider a set of crowd-shippers (suppliers) $s \in S$, which is a realization drawn from random variable $\xi_S$. Thereby, we consider the set of potential hubs $H^{p}$ and binary decision variables $o_h$ for all $h \in H^{p}$ equal to 1 if hub $h$ is opened and 0 otherwise. A hub can be opened at a cost $f$ and a maximum of $Q$ hubs can be opened. The fixed cost of opening a hub mainly consists of the daily rental costs of a location and the maintenance costs of parcel lockers. Acquiring the parcel lockers is a one-time investment and is therefore neglected. The costs involved with the second and third stage depend on these decisions, as well as the realization of demand and supply and are denoted as $\psi_2(o, D)$ and $\psi_3(o, w, D, S)$ respectively. The first stage objective is to minimize the sum of the costs of opening hubs and the expected costs of the second and third stage. We denote with $\mathbb{E}_{\xi}[\cdot]$ is the expected value function over random variable $\xi$. The first stage can be formulated as follows:
\begin{align}
    \min f \sum_{h \in H^{p}} o_h &+ \mathbb{E}_{\xi_D}[\psi_2(o, D)], \label{eq:stage1_1}\\
    s.t. \sum_{h \in H^{p}} o_h &\leq Q,\\
    o_h &\in \mathbb{B} &\forall h \in H^{p}.
\end{align}

\noindent In the second stage we decide which demand unit to assign to which hub. For this, we introduce binary decision variables $w_{dh}$ which is equal to 1 if demand unit $d \in D$ is assigned to hub $h \in H^{p}$, and 0 otherwise. We assume there are no costs involved with the second stage (at least, there is no direct cost difference between assigning to different hubs) other than the expected costs of the third stage. Every demand unit can be assigned to at most one hub and this hub has to be open. This leads to the following formulation of the second stage:
\begin{align}
    \psi_2(o, D)& = \mathbb{E}_{\xi_S}[\psi_3(o, w, D, S),]\\
    s.t. \sum_{h \in H^{p}} w_{dh} &\leq 1 &\forall d \in D,\\
    w_{dh} &\leq o_h &\forall d \in D, h \in H^{p},\\
    w_{dh} &\in \mathbb{B} &\forall d \in D, h \in H^{p}.
\end{align}

\noindent Finally, in the third stage, the parcels are assigned to crowd-shippers. For this, we introduce binary decision variables $x_{ds}$, which is equal to 1 if demand unit $d \in D$ is assigned to crowd-shipper $s \in S$, and 0 otherwise. The crowd-shipper fee for every parcel is denoted as $p_{\tup {cs}}$, and by definition is a fixed cost per parcel. The cost of regular delivery is defined as $p_{\tup {reg}}$ and comprises all costs associated to last-mile delivery such as fuel cost, driver salary and cost of maintenance and repair. Clearly, a parcel can be assigned to only one crowd-shipper and a crowd-shipper can be assigned at most one parcel. Thereby, a parcel can only be assigned to a crowd-shipper if they can feasibly pickup and deliver this parcel, given the hub the parcel was assigned to in the second stage. We use a binary parameter $u_{dsh}$ which is equal to 1 if crowd-shipper $s \in S$ can feasibly pickup parcel $d \in D$ from hub $h \in H^{p}$ and deliver it to the final destination of the parcel. Potential crowd-shippers are assumed to have a maximum detour $\tau$ they are willing to make to pickup and deliver a parcel, which defines this feasibility parameter. The third cost can then be formulated as:
\begin{align}
    \psi_3(o, w, D, S) = p_{\tup {reg}} |D| &+ (p_{\tup {cs}} - p_{\tup {reg}}) \sum_{d \in D} \sum_{s \in S} x_{ds} \label{eq:stage3_1},\\
    \sum_{d \in D} x_{ds} &\leq 1 &\forall s \in S,\\
    \sum_{s \in S} x_{ds} &\leq 1 &\forall d \in D,\\
    x_{ds} &\leq \sum_{h \in H^{p}} w_{dh}u_{dsh} &\forall d\in D, s \in S,\\
    x_{ds} &\in \mathbb{B} &\forall d\in D, s \in S. \label{eq:stage3_last}
\end{align}

\noindent The difficulties of solving this problem in \eqref{eq:stage1_1} - \eqref{eq:stage3_last} are three-fold. First, we are dealing with two separate layers of uncertainty. The first layer of uncertainty is in the demand; with next-day delivery being extremely common, the number of demand requests in every region is uncertain up to a day before delivery. The second layer of uncertainty is the supply. The number of crowd-shippers and their itineraries are generally uncertain up to shortly before departure of the crowd-shipper. Whereas some crowd-shippers may know there schedule well in advance, others may only make themselves available a few minutes before departure. Thereby, exact schedules may be prone to last-minute changes. The second difficulty is that decisions in each of the three stages are heavily intertwined. Decisions in the first and second stage are made based on expected costs and actions in the third stage, whereas the optimal third-stage decisions and corresponding costs depend on the decisions that were made in the first and second stage. This illustrates the importance of solving the problem as a whole and the inability to decompose it. The third difficulty is that the large size of the problem in urban areas causes a computational burden.
\\\\
These difficulties combined make it impossible for the problem to be solved exactly for a realistic case-study. Therefore, in this paper, we approximate the third-stage operational cost using a continuum approximation approach. The second-stage is relaxed and we assume that a parcel can be picked up from any hub. Later, the distribution of parcels over hubs can be determined to be compliant with the first and third-stage decisions. Based on the approximated third-stage costs, we optimize the first stage strategic decisions. By using an approximation of the second and third stage, we are able to evaluate many potential hub location combinations within a reasonable amount of time and therefore explore a large search space.

\subsection{Continuum Approximation}\label{subsec:ca}
\begin{figure}[H]
    \centering
    \includegraphics[width=0.65\textwidth]{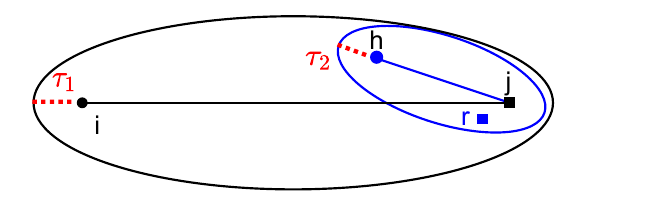}
    \caption{Geometric interpretation of matching problem}
    \label{fig:geometric}
\end{figure}

To approximate the third stage costs, as well as estimate the demand served in every region of the network, we use a CA approach. We consider a network split into $R$ regions. Daily demand for small parcels in every region $r \in R$ is equal to $d_r$. Potential crowd-shippers travel between regions such that the daily average number of crowd-shippers with origin $i$ and destination $j$ is equal to $\lambda_{ij}$, for all $i,j \in R$. The distance between those regions is denoted as $t_{ij}$. Potential crowd-shippers are assumed to have a maximum detour $\tau$ they are willing to make to pickup and deliver a parcel, in return for a reward $p_{\tup {cs}}$. We note that in reality, the reward may depend either on the distance for which the parcel is carried, or the exact detour made by the crowd-shipper and $\tau$ may be heterogeneous among all crowd-shippers. Whenever a parcel is not delivered by a crowd-shipper, it has to be delivered by a regular delivery van at a cost $p_{\tup {reg}}$. We define parameter $e_{ijhr}$ which is equal to 1 if a crowd-shipper with origin at $i$ and destination at $j$ can pickup a parcel at hub $h$ and deliver it to the final destination in region $r$, and 0 otherwise. When $e_{ijhr} = 1$, this is referred to as a feasible assignment. We remark the relation between the parameter $e_{ijhr}$ and the parameter $u_{dsh}$ in the discrete formulation. One of the main advantages of the continuum approximation approach is that the computation of the variables $e_{ijhr}$ relies only on the size of the network (that is, the number of regions) and not on the number of crowd-shippers nor the number of parcels. Thereby, it is independent of the realizations of demand and supply, whereas $u_{dsh}$ depends on the realizations of $D$ and $S$.
\\\\
Geometrically, the matching problem and the definition of $e_{ijhr}$ can be interpreted through two ellipses, as depicted in Figure \ref{fig:geometric}. The original route from origin to destination of a crowd-shipper is depicted by the black line. The crowd-shipper has to make a detour by performing a pickup at $i$ and a delivery at the final destination of the parcel $r$. For this to be feasible within the maximum detour $\tau$, the following should hold. The hub location $h$ should lie within the ellipse with focus points $i$ and $j$, where the distance between the focus and the closest vertex is equal to $\tau_1$. Thereby, the final destination $r$ should lie within a second ellipse with focus points $h$ and $j$, where the distance between the focus and the closest vertex is equal to $\tau_2$. If $\tau_1 + \tau_2 \leq \tau$, the pickup and delivery can be made within the maximum detour. The difficulty of this problem is marked by the variability in the values of $\tau_1$ and $\tau_2$ under the constraint $\tau_1 + \tau_2 \leq \tau$. In addition to this, there is a dependency between the size of the two ellipses. As the point $h$ has to lie within the black ellipse, it restricts the size of the blue ellipse. We also note that we cannot omit the second ellipse by simply considering the line segment $|hj|$ inside the black ellipse in stead of just point $h$, because of the influence of the direction of this segment on the total detour. Due to these difficulties, we resort to numerical approximations.

\FloatBarrier
\FloatBarrier

\noindent When multiple hubs are opened, the main difficulty is that both demand and supply have to be split over the various hubs. Only when the expected number of parcels delivered by a specific hub is independent of the other hubs, this problem can be separated into subproblems. However, especially when two potential hubs are close together, it is clear that this independency is not true. When multiple hubs are opened, each parcel can be served and each crowd-shipper can pick up parcels from multiple hubs as long as they are within their maximum detour $\tau$. We define $\Tilde{e}_{ijr}$ which is equal to 1 if a crowd-shipper with origin at $i$ and destination at $j$ can pickup a parcel \textit{from at least one open hub} and deliver it to the final destination in region $r$, and 0 otherwise. Specifically:
\begin{equation}
    \Tilde{e}_{ijr} = \min(1, \sum_{h \in H}e_{ijhr}).
\end{equation}
We want to approximate the number of parcels delivered from this hub to all regions separately, as well as the total number of parcels that can be delivered through crowd-shipping. To approximate this, we let $\hat{d}_r$ be the expected demand in region $r$ and let $\hat{\lambda}_{ij}$ the expected supply from $i$ to $j$. Then, we define $s_{ij} = \sum_{r \in R}\Tilde{e}_{ijr}\hat{d}_r$ as the total demand that can potentially be served by crowd-shippers going from $i$ to $j$. For the sake of the approximation, we assume that a crowd-shipper is equally likely to choose any of the parcels they can feasibly deliver. Following from this, the probability that he picks a parcel with destination region $r$ is equal to $\frac{\hat{d}_r}{s_{ij}}$ if $\Tilde{e}_{ijr} = 1$ and 0 otherwise.  We can then consider all potential crowd-shippers to obtain the following estimated served demand: 
\begin{align}
    y_r &= \sum_{i,j \in R} \Tilde{e}_{ijr} \hat{\lambda}_{ij} \frac{\hat{d}_r}{s_{ij}} &\forall r \in R. \label{eq:ca_y}
\end{align}
Here, we sum over all origin-destination regions of potential crowd-shippers. Note that this is not necessarily the same as the origin-destination regions of the parcels. An origin-destination region is only considered if the assignment is feasible. Every crowd-shipper in $\hat{\lambda}_{ij}$ is considered to make the delivery with the same probability $\frac{\hat{d}_r}{s_{ij}}$. If no crowd-shippers can be assigned to a region $r$, that is either $\hat{\lambda}_{ij} = 0$ or $\Tilde{e}_{ijr}=0$ for every origin-destination pair, then the number of delivered parcels will always be zero. Similarly, if the expected demand $\hat{d}_r$ is zero, $y_r$ will also be zero. 
\\\\
If $\hat{d}_r$ is non-zero, it is possible that crowd-shippers with different origin-destination pairs are assigned to the same parcel-destination region $r$. As this could lead to an overestimation of served demand in that region ($y_r > \hat{d}_r$), we take into account that at most $\hat{d}_r$ demand can be delivered to a region $r$. Therefore, we set the estimated number of parcels delivered to
\begin{align}
    z_r &= \min(\hat{d}_r, y_r) &\forall r \in R.
\end{align}
Especially if supply is high, by overestimating $y_r$ in region $r$ (i.e. $y_r > \hat{d}_r$), it is likely that $y_{r'}$ for another region $r' \neq r$ will be underestimated. Therefore, we use an iterative process to ensure that this unavoidable overestimation is accounted for in the other regions. We consider the leftover demand $l_r = \max(0, y_r - \hat{d}_r)$ and split it evenly over the potential suppliers. Similar to the assignment of parcels to crowd-shippers, we assume that every crowd-shipper that can be feasibly assigned to a region $r$ (that is, every crowd-shipper with origin $i$ and destination $j$ for which $\Tilde{e}_{ijr} = 1$), is equally likely to be assigned to one of the leftover demand units in $l_r$. Therefore, the $l_r$ leftover demand units are split over the origin-destination pairs proportional to the number of suppliers that could be feasibly assigned to region $r$. We define the leftover supply as follows:
\begin{equation}
    \hat{\lambda}'_{ij} = \sum_{r \in R} l_r\frac{\Tilde{e}_{ijr}\hat{\lambda}_{ij}}{\sum_{i,j \in R} \Tilde{e}_{ijr}\hat{\lambda}_{ij}}. \label{eq:ca_lambdaprime}
\end{equation}
Thereby, we define the unserved demand $\hat{d}'_r = \hat{d}_r - z_r$. All demand units that are already expected to be served by previously assigned crowd-shippers no longer need to be considered and are therefore disregarded. We then compute $y_r$ according to Equation \eqref{eq:ca_y}, but now using  $\hat{\lambda}'_{ij}$ and $\hat{d}'_r$ as inputs in stead of $\hat{\lambda}_{ij}$ and $\hat{d}_r$. Using these we find an additional portion of demand which can be served and we update the estimated demand served and the leftover demand as follows: 
\begin{align}
    z_r &= \min(\hat{d}_r, z_r + y_r) &\forall r \in R.\\
    l_r &= \max(0, y_r - \hat{d}'_r) &\forall r \in R.
\end{align}
This iterative process can be repeated until the leftover demand $l_r$ is zero for all regions $r \in R$. Intuitively, if $l_r = 0$, all previously overestimated demand has been compensated for. However, we emphasize that the simplifying assumptions of proportional assignment in Equation \eqref{eq:ca_y} and \eqref{eq:ca_lambdaprime} can lead to suboptimal assignments. If $l_r = 0$ for all regions already after the first iteration, the estimated served demand over all regions reduces to $\sum_{i,j \in R} \hat{\lambda}_{ij}$, as is shown in Equation \eqref{eq:optimal_derivation}. This is typically the case if supply is much lower than demand, such that served demand is never overestimated, and when asymmetries are relatively low. In this case, it is clear that our approximation is equal to the total available supply. As not more parcels can be delivered by crowd-shippers than there are crowd-shippers available, our approximation of the number of parcels delivered by crowd-shippers is equal to the maximum (optimal) number of parcels delivered by crowd-shippers. Another extreme case exists when supply is extremely large. In this case, the approximation will lead to $z_r = \hat{d}_r$, which corresponds to all demand being served. 
\begin{align}
    \sum_{r \in R} z_r &= \sum_{r \in R} y_r = \sum_{r \in R} \sum_{i,j \in R} \Tilde{e}_{ijr} \hat{\lambda}_{ij} \frac{\hat{d}_r}{s_{ij}} \nonumber \\
    &= \sum_{r \in R} \sum_{i,j \in R} \Tilde{e}_{ijr} \hat{\lambda}_{ij} \frac{\hat{d}_r}{\sum_{r \in R}\Tilde{e}_{ijr}\hat{d}_r} = \sum_{i,j \in R} \hat{\lambda}_{ij} \frac{ \sum_{r \in R} \Tilde{e}_{ijr} \hat{d}_r}{\sum_{r \in R}\Tilde{e}_{ijr}\hat{d}_r} \label{eq:optimal_derivation} \\
    &= \sum_{i,j \in R} \hat{\lambda}_{ij}. \nonumber 
\end{align}

\noindent The total cost can be obtained directly from the results of the approximation. Let $z_r(H)$ be the number of parcels delivered to region $r$ if the selection of hubs $H$ is opened. The total approximated cost, similar to that defined in Section \ref{subsec:discrete}, is as follows:
\begin{equation}
    C_{\tup {total}}(H) = f|H| + p_{\tup {cs}} \sum_{r \in R} z_r(H) + p_{\tup {reg}} \sum_{r \in R} (d_r - z_r(H)).
\end{equation}

\subsection{Stage 1: Determining Hub Locations}\label{subsec:stage_1}
The cost for a given set of hub locations can be determined efficiently using the approximation methods discussed in the previous section. Despite this, in highly populated urban networks the number of options to consider can still be extremely large. Let $n$ be the number of possible hub locations and let $Q$ be the maximum number of hubs that can be built. Enumerating all possible hub locations will have complexity $\mathcal{O}(\frac{n!}{Q!(n - Q)!})$ when exactly $Q$ hubs have to be opened and even larger when at most $Q$ hubs have to be opened. If there is no limit on the maximal number of hubs, complexity is $\mathcal{O}(2^n)$. Clearly, enumerating all options is impossible for large networks and therefore an efficient heuristic will have to be used to determine the best hub locations. We propose a Large Neighborhood Search (LNS) heuristic to solve this problem. LNS heuristics explore a complex neighborhood to find better candidate solutions \citep{pisinger2010large}. We explore the neighborhood in an efficient way by using metrics for the quality of solutions, that allow to select candidates to \textit{destroy} and \textit{repair} in a smart way.   

\subsubsection*{Initialization} 
To initialize the heuristic, we precompute a number of components. First, we compute the 4-dimensional matrix $E$ with elements $e_{ijhr}$ for every potential hub. Thereby, we compute the single hub objectives for each hub $h \in H^p$ which will be used as a quality metric for the hubs. This metric will be referred to as $v_h$, which is equal to $C_{\tup {total}}(\{h\})$. The intuition behind this is that a hub that performs well on its own is more likely to perform well in combination with other hubs. Nevertheless, hubs that are serving mostly similar demand locations and attract similar crowd-shippers might perform poorly if they operate together, as one hub has little added value over the other hub. For this reason, we construct a similarity measure $s_{h_1,h_2}$ for how similar two hubs are in terms of the service area of crowd-shippers using that hub. Specifically, $s_{h_1,h_2}$ is determined as follows:

\begin{equation}
    s_{h_1,h_2} = \frac{[\sum_{i,j,r \in R} \min(e_{ijh_1r}, e_{ijh_2r}) \lambda_{ij}]^2 } {[\sum_{i,j,r \in R} e_{ijh_1r} \lambda_{ij}][ \sum_{i,j,r \in R} e_{ijh_2r} \lambda_{ij}]}.
\end{equation}

\noindent Clearly, two hubs that are very similar in terms of service area are likely to have a lower gain in performance when they are combined. For this reason, this similarity measure will be used to select dissimilar hubs to be combined. 
\\\\
We use a multi-start heuristic, so that we randomly determine $\eta$ initial solutions. Every initial solution is generated according to a simple construction heuristic. Every hub is randomly chosen with a probability proportional to the quality of the hub in the single-hub solution. By using a multi-start heuristic we aim to increase the search space and therefore decrease the likelihood of ending up at a local optimum. 

\subsubsection*{Body of algorithm} 
We use a large neighborhood search algorithm for a fixed number of $\kappa$ iterations. In every iteration, we use a number of operations on the current solution $\Omega$ and obtain the corresponding objective value. A newly generated solution is always accepted if it is an improvement over the previous solution. We consider the following operators: 
\begin{enumerate}
    \item \textit{$O_1$ - Repair operator:} For every hub $h$ that is not in the current solution $\Omega$, we compute the following metric: $\frac{v_h^\alpha}{[\sum_{\omega \in \Omega} s_{h, \omega}]^\beta}$ where $\alpha$ and $\beta$ are tuning parameters that determine the relative importance of single-hub performance and inter-hub similarity. This metric determines the best hub to be added to the current solution, taking into account the quality of the hub in the single-hub solution as well as the similarity to the other hubs. Here, we use the sum of the similarity with all hubs in the current solution $\Omega$. The intuition behind this is that a hub that is similar to two hubs rather than only one should have a lower metric value. Alternatively, the maximum similarity across all hubs in the current solution $\Omega$ can be used in the denominator to replace the sum. A new hub is added randomly with a probability proportional to the value of this metric. 
    \item \textit{$O_2$ - Destroy operator:} For every hub $h$ that is in the current solution $\Omega$, we determine the same metric as in $O_1$ and randomly drop a hub with a probability proportional to the inverse of the metric in $O_1$. By taking the inverse, hubs that are very similar to other hubs as well as hubs that have relatively low single-hub performance are most likely to be removed. 
    \item \textit{$O_3$ - Swapping operator:} A sequential combination of $O_1$ and $O_2$ where a hub in the current solution $\Omega$ is replaced by another hub that is not in the current solution. We first destroy a hub $\omega \in \Omega$ using destory operator $O_2$ and then add a hub using repair operator $O_1$, based on the similarity with the remaining hubs $\Omega\setminus\{\omega\}$.
\end{enumerate}

\noindent We note the high level of randomness in our algorithm. The reason for this is to stimulate diversification. Generally, using a heuristic to determine hub locations rather than directly solving an optimization problem leads to more robust solutions. Given the large level of uncertainty in our problem, such a robust solution is of importance. 

\subsection{Stage 2: Parcel-Hub Assignment}\label{subsec:stage_2}
After the hubs are determined in stage 1, parcels have to be assigned to hubs on a day-to-day basis. At this stage, demand is assumed to be known exactly (all orders of parcels have been collected), but the exact supply realizations are still unknown (crowd-shippers can announce their availability last-minute). By simply assigning parcels to the closest hub in terms of distance, the importance of potential supply-flow is neglected. For the sake of comparison, we consider a distance-based metric that assigns all the demand of a region to the opened hub that is closest to that region. Consider the set of opened hubs $H$ and consider $d_r$ demand units with destinations in region $R$. We define $a_{rh}$ the number of parcels with final destination $r \in R$ that are stored at hub $h \in H$ as follows:

\begin{equation}
    a_{rh} = d_r \mathbbm{1}_{\Big[t_{rh} = \min\limits_{h' \in H} t_{rh'}\Big]}.
\end{equation}

\noindent Although a hub can be close in terms of distance, if very few crowd-shippers are able to feasibly deliver a parcel from that hub to the final destination, such a parcel-hub assignment can perform poorly. Therefore, we develop an assignment strategy based on the CA-estimates obtained using the algorithm as described in Section \ref{subsec:ca}.   We solve the single-hub-approximation for every hub $h \in H$ to obtain the expected number of parcels delivered to region $r \in R$, $z_r^h$, where the superscript is now added to distinguish between hubs. We then assign the parcels proportional to the expected number of parcels delivered from a specific hub location. In this case, we define $a_{rh}$ as follows:
\begin{equation}
    a_{rh} = d_r\frac{(z_r^h)^\gamma}{\sum_{h \in H} (z_r^h)^\gamma}.
\end{equation}
\noindent The tuning parameter $\gamma$ can be used to give extra weight to larger hubs (hubs with high expected demand) and less weight to smaller hubs (hubs with lower expected demand).  

\subsection{Stage 3: Parcel-Crowd-Shipper Matching}\label{subsec:stage_3}
In the third stage of our problem, when the hub locations are known and parcels are distributed over these hub locations, parcels have to be assigned to crowd-shippers. Generally, parcels are matched to crowd-shippers dynamically, upon arrival of the crowd-shippers. A parcel can only be matched to a crowd-shipper if they can pickup and deliver the parcel within their maximum detour, this is referred to as a \textit{feasible match}. If a crowd-shipper can be feasibly matched to multiple parcels, the operator has to decide which parcel to assign to the crowd-shipper in order to maximize the total number of delivered parcels over the entire planning horizon. This decision is made under uncertainty of future arriving crowd-shippers. We consider the following four alternative matching approaches. The static matching approach and the batch matching approach rely on solving an integer linear programming problem, whereas the minimal-detour and CA-based matching use a simple decision rule. Furthermore, the static and batch matching approach use information about future crowd-shippers, whereas the other two approach do not integrate any predictions.  

\subsubsection*{Static Matching}
The static matching approach assumes complete knowledge of all future crowd-shippers that will arrive. The static matching can be obtained by solving an ILP problem, that can be taken from the third stage of the stochastic programming formulation \eqref{eq:stage3_1} - \eqref{eq:stage3_last}. Here, the hubs are fixed and the demand-hub assignment has been made. Therefore, we can eliminate the variables $w_{dh}$ and directly incorporate the hub as the origin for a demand. Let $D$ be the set of demand requests and $S$ the set of potential suppliers (i.e. crowd-shippers). Decision variable $x_{ds}$ is equal to 1 if demand unit $d \in D$ is matched to supplier $s \in S$ and 0 otherwise. Parameter $u_{ds}$ is equal to 1 if $(d,s)$ is a feasible match and 0 otherwise. As the assigned hub is fixed, the subscript $h$ is omitted here. The static matching problem can then be formulated as follows:

\begin{align}
    \max \sum_{d \in D}\sum_{s \in S} x_{ds}, \label{eq:static_objective}\\
    \sum_{d \in D} x_{ds} &\leq 1 &\forall s \in S, \label{eq:static_supply}\\
    \sum_{s \in S} x_{ds} &\leq 1 &\forall d \in D, \label{eq:static_demand}\\
    x_{ds} &\leq u_{ds} &\forall d\in D, s\in S, \label{eq:static_feasible}\\
    x_{ds} &\in \mathbb{B}  &\forall d\in D, s\in S. \label{eq:static_range}
\end{align}

\noindent The objective in \eqref{eq:static_objective} is to maximize the total number of matched parcels. For this, we note that the cost parameters can be disregarded, because the hubs are fixed. Constraints \eqref{eq:static_supply} and \eqref{eq:static_demand} enforce that every crowd-shipper can be matched to at most one parcel and every parcel can be matched to at most one crowd-shipper, respectively. Constraints \eqref{eq:static_feasible} ensure that all matches are feasible and Equation \eqref{eq:static_range} defines the binary range of the decision variables. 
\\\\
We emphasize that, despite the problem being of such a nature that solving the Linear Programming (LP) approximation of the problem gives the optimal integer solution, building and solving the model can still take considerable time. The main reason for this is the sheer size of the sets $D$ and $S$ that make the problem expensive to build, especially through the parameters $u_{ds}$. 
\\\\
As the static matching is made without uncertainty about the future crowd-shippers, this forms an upperbound to the other matching approaches. In reality, this could be very unrealistic if crowd-shippers make themselves known only shortly before departing. Despite this, it forms a useful benchmark to compare the dynamic assignment strategies and the continuum approximation. 

\subsubsection*{Batch Matching}
Batch matching is a popular matching approach in the sharing economy literature. For example, \cite{syed2019asynchronous} use batch optimization to match passengers and drivers in ride-hailing services. Batch matching refers to a matching approach where a group (i.e. batch) of crowd-shippers is matched at the same time. The advantage of batch matching is that information on more than one crowd-shipper can be gathered to improve the quality of the match. The disadvantage is that it requires multiple potential crowd-shippers to have indicated their availability and therefore it requires potential crowd-shippers to announce their participation well before they leave for their trip. In reality, participation can be a last-minute decision for many potential crowd-shippers. Although batch matching is much more realistic than the completely static matching, it is more limiting than the minimal-detour matching and CA-based matching that do not have this requirement. The batch matching is solved using a similar ILP problem as for the static matching, however only a limited number of potential crowd-shippers are used. We use a batch $B$ to replace the total set of suppliers $S$. Thereby, the demand set $D$ only contains parcels that are not yet picked up by or assigned to other crowd-shippers.
\\\\
As we are using a dynamic assignment strategy, computation time is of major importance. We therefore note that the size of the batch $B$ should be small enough for the matching to be found in real time. In addition to this, the time difference between the first and last departing crowd-shipper in the batch should be small enough, such that the first crowd-shipper can wait for his matched parcel to be announced. On the other hand, increasing the batch size generally leads to better service levels, as more information on supply can be used. These components drive the trade-off for the choice for a proper batch size. 


\subsubsection*{Minimal-Detour Matching}
The minimal detour matching matches the crowd-shipper to the parcel for which the crowd-shipper minimizes the detour. This matching strategy therefore does not use any knowledge of future crowd-shippers. This strategy minimizes crowd-shipper inconvenience, but is likely to be suboptimal as it does not use any information on future supply. A minimal-detour matching is a realistic representation of what the matching would look like if crowd-shippers themselves could choose which parcel they want to deliver, if they were given all parcel options. Given that the reward is constant, a crowd-shipper will then always minimize his/her inconvenience. 
\\\\
If we consider that the set $D$ only contains parcels that are not yet assigned to other crowd-shippers, similar to the batch matching, and we consider $s \in S$ the current crowd-shipper. We let $or(\cdot)$ and $des(\cdot)$ be the origin and destination location, respectively, for a crowd-shipper or parcel. Then the minimal detour matching chooses the parcel as follows:
\begin{equation}
    d_{min} = \arg\min_{d \in D} t_{or(s)or(d)} + t_{or(d)dest(d)} + t_{dest(d)dest(s)} - t_{or(s)dest(s)}.
\end{equation}

\subsubsection*{CA-based Matching}
The minimal-detour matching as well as the batch matching only use information on the supplier(s) that is (are) currently available. In this way, a lot of information about potential crowd-shippers that arrive in the future remains unused, although this can be very useful. Generally, historic data is available that provides insights in the expected number of crowd-shippers. The CA-based matching approach aims to exploit the information on the expected number of delivered parcels obtained through the CA to improve the matching quality. Here, an arriving crowd-shipper is assigned the parcel with the demand region that has the lowest expected number of delivered parcels relative to the total demand in that region ($\frac{z_r}{d_r}$). The intuition behind this assignment strategy is that we favor parcels that are less likely to be delivered in the future. By doing so, we increase the total expected number of parcels delivered over the entire planning horizon. 
Compared to batch matching, the CA-based matching only uses a simple metric to determine the assignment, rather than solving an ILP problem. Therefore, the match can be determined extremely fast, similar to the minimal-detour matching, which forms a major advantage for our real-time application. Here, we do not take into account fluctuations of the arrival rates of crowd-shippers during the day. This is marked as an interesting direction of future research.

\section{Results}\label{sec:results}
In this section we evaluate the performance of the developed CA approach to find the optimal hub locations, as well as the potential of hub-based crowd-shipping. Our results are obtained through a case study based on the city of Washington DC, of which the details are described in Section \ref{subsec:case_study}. We develop a discrete event simulator based on this case study, which is described in Section \ref{subsec:simulator}. In Section \ref{subsec:comparison} we compare our CA approach to an exact approach, which forms a performance benchmark. In Section \ref{subsec:comparison_approach} we further evaluate the performance of our CA approach by comparing our solution algorithm to a well-known simulation-optimization approach. In Section \ref{subsec:results_network} we evaluate the results on the network and in Section \ref{subsec:optimal_number_hubs} we perform a sensitivity analysis on the optimal number of hubs. We compare our CA approach to an approach that is non-predictive in all three stages in Section \ref{subsec:parcel_hub_results} and we compare the various parcel-crowd-shipper assignment policies in Section \ref{subsec:assignment_policies}. 

\subsection{Case Study}\label{subsec:case_study}
We use the city of Washington DC and the surrounding metropolitan area as a case study. The city of Washington DC has around 700,000 inhabitants and the entire agglomeration has around 7 million inhabitants. Washington DC hosts one one the biggest Bike-sharing platforms in the USA: Capital Bikeshare. Capital Bikeshare has over 500 stations and 4500 bikes. Such a bike-sharing platform with a large number of users forms a good base for a crowd-shipping service. We consider the locations of the bike-sharing stations as potential hub locations and use them as approximations for regions making up the entire service area. 

\begin{figure}[ht]
    \centering
    \includegraphics[width = 0.55\textwidth]{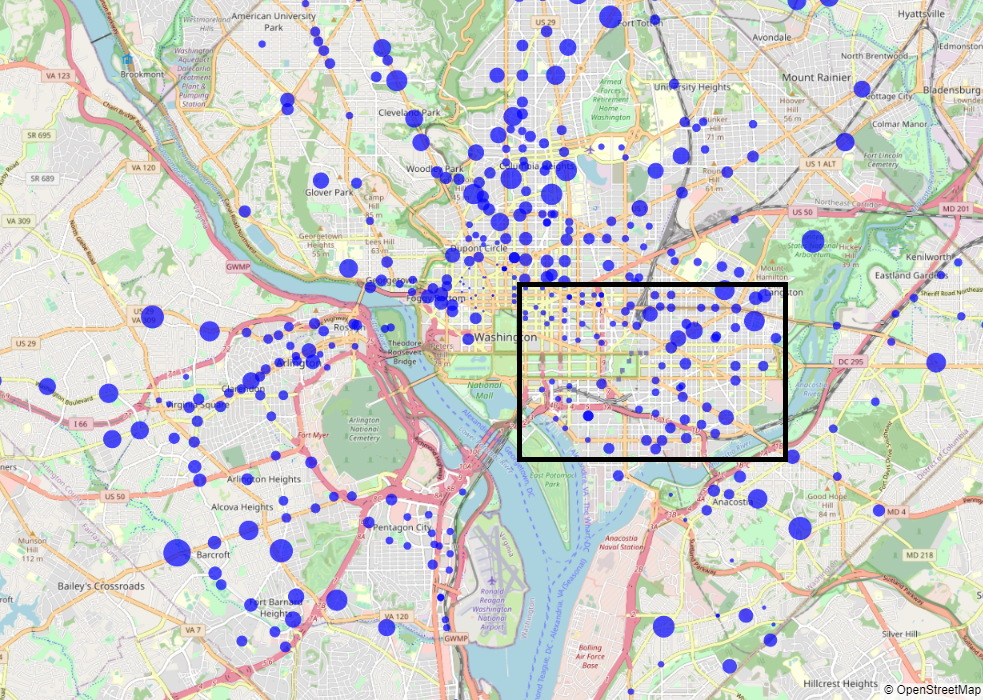}
    \caption{Bubble chart of bike-sharing stations, where the size of the bubble is determined by the population in the area. The area considered in the case study is captured by a black box (5.5$\times$3.5km). }
    \label{fig:stations}
\end{figure}

Based on the station names, approximate coordinates of the locations are extracted from \cite{google2020maps} as are the distances between every pair of stations. In this way, we use a realistic distance, incorporating asymmetries caused by one-way streets. An approximation of the surrounding population has been made using \cite{census2021} data, which in turn has been used as a proxy for demand. Historic system data from the \cite{capital2020} database has been used to identify origin-destination pairs for the supply of crowd-shippers. The stations that are used as regions are displayed in Figure \ref{fig:stations}. This figure displays a bubble chart of those regions, where the size of the bubble is determined relative to the population around the corresponding station. Whereas the most used stations in terms of origins and destination of crowd-shippers are around Union Station, the Mall and the center of Washington DC, demand is higher in the suburbs. This shows the large asymmetry in supply and demand locations that is usually present in crowd-shipping systems in real urban areas, making this case-study especially realistic and interesting. To perform an extensive sensitivity analysis within a reasonable amount of time, we use a selection of 90 stations in the urban area of Washington DC, marked by the black box in Figure \ref{fig:stations}. For the sake of computation time, we only use crowd-shippers that have their origin and destination within the selected area. Any potential crowd-shipper with either an origin or destination outside the selected black box is neglected. In our experiments, demand has been fixed at the expected value such that supply is the only uncertain variable in the problem.

\noindent We use the following parameter values to get a realistic interpretation of the results. We assume that the cost of opening and operating a hub location $(f)$ is equal to 250\$ per day. This is based on the average rental price in Washington DC to place the storage lockers and the maintenance cost involved in operating the hubs. The cost of regular delivery is set to be 15\$, similar to \cite{le2019supply}. A part of these costs is still made to serve the hubs, therefore, we consider a slightly lower value for $p_{\tup {reg}}$ equal to 7.5\$.  Our baseline scenario assumes that crowd-shippers are willing to make a detour ($\tau$) of at most 500 meters and wish to receive ($p_{\tup {cs}}$) 5\$ to make a delivery. According to \cite{statista2021}, the number of online transactions per online shopper in North America is equal to 19 per year. Given the approximated population in the selected region, this is amounts to a daily demand of roughly 4300 parcels per day.
\\\\
The LNS parameters are tuned in order to find a good objective value within reasonable time. The algorithm uses 5 initial solutions ($\eta$) and 500 iterations ($\kappa$). The repair and destroy operators take parameters $\alpha$ equal to 4.5 and $\beta$ equal to 8. These values were chosen based on a grid search over a large range of parameter values. The parameters were chosen such that they maximize the number of multi-start initial solutions for which the best objective was found and at the same time minimize the number of iterations needed to find this objective. The tuning parameter $\gamma$, used for the parcel-hub assignment, is chosen equal to 1.

\subsection{Discrete Event Simulator} \label{subsec:simulator}
We develop a discrete event simulator to simulate the dynamic operational process based on the described methods and case study. The simulator is initialized by generating a set of demand requests consisting of only a destination region and a set of potential suppliers consisting of an origin and destination region and a starting time of the trip. Each simulation consists of 1 day. For the sake of this simulation, crowd-shippers are assumed to make themselves available at the start of their trip. All generated demand requests are assigned to a hub based on one of the strategies described in Section \ref{subsec:stage_2}. All potential crowd-shippers are sorted in ascending order of their start times. The generation of demand and supply is performed using a pseudo-random number generator, such that simulations using various policies can be directly compared.
\\\\
Upon arrival of a crowd-shipper, a parcel is assigned to this crowd-shipper based on one of the strategies described in Section \ref{subsec:stage_3}. For the static assignment strategy, the assignment is pre-defined at initialization, as this strategy is made with full knowledge of supply. For the batch assignment strategy, the assignment of a batch of crowd-shippers is made at the arrival time of the first crowd-shipper of that batch. A batch size of $|B| = 50$ is used. For the minimal detour and CA-based strategy, the assignment is made for every crowd-shipper individually at their time of arrival. After the assignment, the parcel is reserved for the crowd-shipper for pickup and the crowd-shipper departs from his origin to the origin location of the parcel. A new pickup event is scheduled, taking into account the travel time between the origin of the crowd-shipper and the origin of the parcel. As soon as a parcel is assigned to a crowd-shipper, it is no longer available to be assigned to other crowd-shippers, even when it is not yet picked up. 
\\\\
For every pickup event, a delivery event is scheduled taking into account the travel time between the origin and destination of a parcel. As a crowd-shipper can deliver at most one parcel, a crowd-shipper is discarded after the delivery event is completed. For every delivery event, the number of served parcels and total costs are updated and the detour made by the crowd-shipper is stored. The simulation ends when all parcels have been delivered or when all crowd-shippers have either completed a delivery or have failed to be assigned to a parcel.

\subsection{Comparison of Continuum Approximation to Static and Dynamic Assignment Strategies}\label{subsec:comparison}
\begin{figure}[H]
    \centering
    \includegraphics[width = \textwidth]{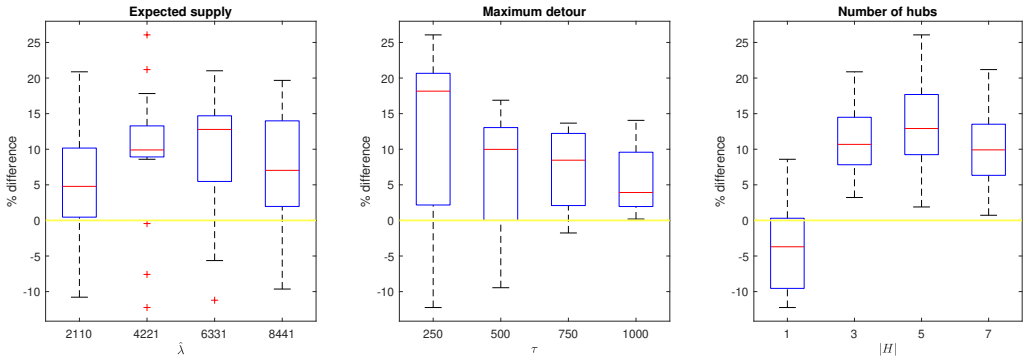}
    \caption{Percentage difference between the predicted served parcels by CA approach and the actual served parcels simulated by the static assignment strategy}
    \label{fig:boxplot_static}
\end{figure}

\begin{figure}[H]
    \centering
    \includegraphics[width = \textwidth]{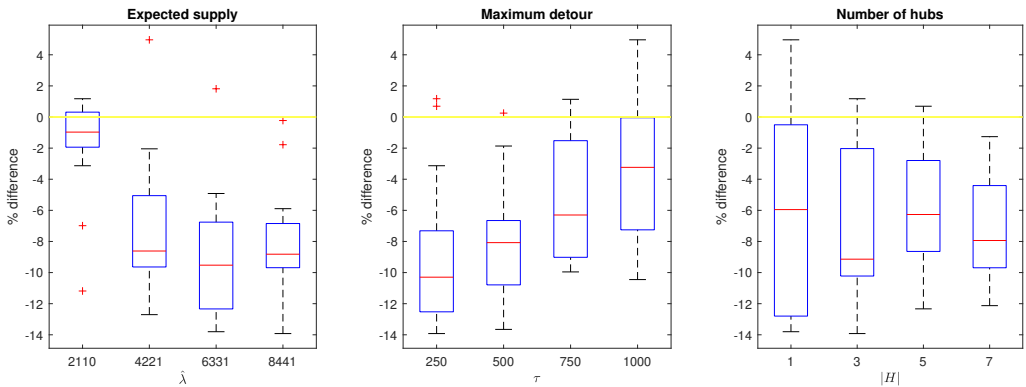}
    \caption{Percentage difference between the predicted served parcels by CA approach and the actual served parcels simulated by the dynamic assignment strategy}
    \label{fig:boxplot_dynamic}
\end{figure}

To evaluate the quality of the CA method, we compare its estimates to two benchmarks. The first benchmark assumes full knowledge of supply \textit{before} parcels are distributed to the hubs.  This means that we use the formulation of the second and third-stage as given by Equations (4) - (12), but for a single realization of demand and supply which is known with certainty. This therefore eliminates the expected value in Equation (4) and allows us to solve the problem to optimality using a standard MIP solver, such as CPLEX. The second benchmark is a realistic dynamic assignment procedure. We use the CA-based parcel-hub assignment and the CA-based matching of parcels and crowd-shippers, which outperforms other assignment strategies (see Sections \ref{subsec:parcel_hub_results} and \ref{subsec:assignment_policies}). A large sensitivity analysis is performed with respect to the expected number potential crowd-shippers ($\hat \lambda$), the maximum detour ($\tau$) and the number of hubs $(|H|)$. The results are given in Table \ref{tab:performance_ca}, which reports the percentage of demand served by crowd-shippers and the percentage deviation from the continuum approximation in parenthesis. The results are visualized in Figures \ref{fig:boxplot_static} and \ref{fig:boxplot_dynamic}, that display the percentage difference between the predicted number of served parcels by CA approach and the actual number of served parcels simulated using the static and dynamic assignment strategy respectively. 

\begin{table}[H]
\caption{Performance of CA compared to  static and dynamic benchmark}
\tiny
\label{tab:performance_ca}
\begin{threeparttable}
\begin{tabular}{ccc|c|rr|rr} \toprule 
$\hat \lambda$ & $\tau$ & $|H|$ & CA & \multicolumn{2}{c}{Static benchmark } & \multicolumn{2}{c}{Dynamic benchmark} \\ \midrule
2110   & 250  & 1    & 12.1 & 10.8         & (-10.8)         & 10.8          & (-11.2)         \\
2110   & 250  & 3    & 23.2 & 28.0         & (20.9)          & 23.4          & (1.2)           \\
2110   & 250  & 5    & 29.9 & 36.1         & (20.5)          & 30.2          & (0.7)           \\
2110   & 250  & 7    & 35.6 & 40.6         & (14.0)          & 34.5          & (-3.1)          \\
2110   & 500  & 1    & 15.7 & 14.6         & (-7.0)          & 14.6          & (-7.0)          \\
2110   & 500  & 3    & 30.5 & 34.0         & (11.3)          & 29.9          & (-2.0)          \\
2110   & 500  & 5    & 38.5 & 42.0         & (9.1)           & 38.6          & (0.3)           \\
2110   & 500  & 7    & 42.0 & 44.4         & (5.6)           & 41.2          & (-1.9)          \\
2110   & 750  & 1    & 18.9 & 18.8         & (-0.2)          & 18.7          & (-0.8)          \\
2110   & 750  & 3    & 35.6 & 37.8         & (6.2)           & 36.0          & (1.1)           \\
2110   & 750  & 5    & 41.9 & 44.5         & (6.2)           & 41.4          & (-1.2)          \\
2110   & 750  & 7    & 46.0 & 47.1         & (2.3)           & 45.5          & (-1.3)          \\
2110   & 1000 & 1    & 22.7 & 22.8         & (0.2)           & 22.8          & (0.1)           \\
2110   & 1000 & 3    & 39.9 & 41.5         & (3.9)           & 40.1          & (0.4)           \\
2110   & 1000 & 5    & 46.6 & 47.5         & (1.9)           & 46.3          & (-0.7)          \\
2110   & 1000 & 7    & 48.2 & 48.5         & (0.7)           & 47.5          & (-1.3)          \\
4221   & 250  & 1    & 23.4 & 20.5         & (-12.2)         & 20.4          & (-12.7)         \\
4221   & 250  & 3    & 44.4 & 52.3         & (17.8)          & 40.2          & (-9.4)          \\
4221   & 250  & 5    & 56.6 & 71.4         & (26.1)          & 51.5          & (-9.1)          \\
4221   & 250  & 7    & 62.6 & 75.8         & (21.2)          & 57.4          & (-8.3)          \\
4221   & 500  & 1    & 30.1 & 27.8         & (-7.6)          & 27.6          & (-8.2)          \\
4221   & 500  & 3    & 55.2 & 60.7         & (10.0)          & 49.4          & (-10.6)         \\
4221   & 500  & 5    & 65.0 & 73.6         & (13.3)          & 57.8          & (-11.0)         \\
4221   & 500  & 7    & 71.4 & 78.5         & (10.0)          & 64.3          & (-9.9)          \\
4221   & 750  & 1    & 36.2 & 36.0         & (-0.4)          & 35.4          & (-2.2)          \\
4221   & 750  & 3    & 62.1 & 67.8         & (9.2)           & 56.5          & (-9.0)          \\
4221   & 750  & 5    & 72.4 & 82.0         & (13.2)          & 68.3          & (-5.7)          \\
4221   & 750  & 7    & 77.7 & 85.3         & (9.8)           & 70.7          & (-9.0)          \\
4221   & 1000 & 1    & 41.4 & 45.0         & (8.6)           & 43.5          & (5.0)           \\
4221   & 1000 & 3    & 71.1 & 78.0         & (9.7)           & 69.6          & (-2.0)          \\
4221   & 1000 & 5    & 78.5 & 87.9         & (12.0)          & 75.0          & (-4.4)          \\
4221   & 1000 & 7    & 80.7 & 88.5         & (9.7)           & 76.1          & (-5.7)          \\
6331   & 250  & 1    & 33.2 & 29.4         & (-11.2)         & 28.6          & (-13.8)         \\
6331   & 250  & 3    & 59.6 & 68.8         & (15.3)          & 52.3          & (-12.3)         \\
6331   & 250  & 5    & 68.9 & 83.4         & (21.0)          & 60.4          & (-12.3)         \\
6331   & 250  & 7    & 75.5 & 90.3         & (19.5)          & 66.4          & (-12.1)         \\
6331   & 500  & 1    & 41.5 & 39.2         & (-5.6)          & 36.1          & (-13.2)         \\
6331   & 500  & 3    & 71.3 & 78.4         & (10.1)          & 62.0          & (-13.0)         \\
6331   & 500  & 5    & 78.0 & 91.2         & (16.9)          & 72.9          & (-6.6)          \\
6331   & 500  & 7    & 82.7 & 93.5         & (13.0)          & 76.4          & (-7.6)          \\
6331   & 750  & 1    & 49.6 & 48.7         & (-1.8)          & 47.1          & (-4.9)          \\
6331   & 750  & 3    & 76.5 & 87.0         & (13.7)          & 71.3          & (-6.9)          \\
6331   & 750  & 5    & 85.4 & 96.1         & (12.6)          & 78.3          & (-8.3)          \\
6331   & 750  & 7    & 88.0 & 99.5         & (13.0)          & 79.3          & (-10.0)         \\
6331   & 1000 & 1    & 57.8 & 60.1         & (3.9)           & 58.9          & (1.8)           \\
6331   & 1000 & 3    & 84.6 & 96.5         & (14.1)          & 76.9          & (-9.1)          \\
6331   & 1000 & 5    & 91.0 & 99.6         & (9.5)           & 85.7          & (-5.9)          \\
6331   & 1000 & 7    & 93.0 & 99.5         & (7.0)           & 83.2          & (-10.5)         \\
8441   & 250  & 1    & 42.2 & 38.1         & (-9.6)          & 36.8          & (-12.9)         \\
8441   & 250  & 3    & 70.2 & 80.7         & (14.9)          & 60.5          & (-13.9)         \\
8441   & 250  & 5    & 77.0 & 91.2         & (18.5)          & 71.6          & (-7.0)          \\
8441   & 250  & 7    & 80.1 & 95.9         & (19.7)          & 74.0          & (-7.6)          \\
8441   & 500  & 1    & 51.2 & 46.4         & (-9.5)          & 44.2          & (-13.7)         \\
8441   & 500  & 3    & 78.4 & 83.4         & (6.4)           & 71.0          & (-9.3)          \\
8441   & 500  & 5    & 84.3 & 97.5         & (15.6)          & 77.6          & (-7.9)          \\
8441   & 500  & 7    & 87.1 & 98.5         & (13.1)          & 81.3          & (-6.7)          \\
8441   & 750  & 1    & 59.4 & 60.5         & (1.9)           & 58.3          & (-1.8)          \\
8441   & 750  & 3    & 86.8 & 97.1         & (11.9)          & 78.8          & (-9.2)          \\
8441   & 750  & 5    & 90.4 & 98.9         & (9.4)           & 82.2          & (-9.0)          \\
8441   & 750  & 7    & 92.4 & 99.5         & (7.7)           & 83.6          & (-9.5)          \\
8441   & 1000 & 1    & 68.7 & 68.9         & (0.4)           & 68.5          & (-0.2)          \\
8441   & 1000 & 3    & 94.4 & 97.4         & (3.2)           & 85.1          & (-9.9)          \\
8441   & 1000 & 5    & 96.3 & 99.6         & (3.5)           & 90.6          & (-5.9)          \\
8441   & 1000 & 7    & 97.7 & 99.6         & (2.0)           & 89.3          & (-8.6)            \\ \bottomrule
\end{tabular}
 \begin{tablenotes}[para]
\footnotesize{The first three columns give the properties of the instance. The supply $\hat{\lambda}$ gives the number of crowd-shippers, the maximum detour $\tau$ is given in meters and $|H|$ gives the number of hubs opened. The last three columns give the percentage of demand served by crowd-shippers, where the number in parenthesis is the difference between the CA approximation and the static benchmark and dynamic simulation benchmark respectively.}
\end{tablenotes}
\end{threeparttable}

\end{table}

\noindent Due to the full knowledge of demand and supply in the static benchmark, the static assignment method provides an upper bound for the (more realistic) dynamic CA-based assignment strategy. The percentage of demand served as estimated by the continuum approximation method is generally between the two benchmarks. The closeness to the realistic benchmark of the simulated dynamic assignment suggests that our CA method provides an accurate approximation of the actual number of parcels that crowd-shippers can deliver. Also, we emphasize that a more advanced dynamic assignment strategy could likely improve the performance further, bringing it closer to the CA value. 
\\\\
The sensitivity analysis indicates that for higher number of potential crowd-shippers, performance slightly deteriorates. The most plausible reason for this is that as the number of potential crowd-shippers increases, the number of possible matches for every parcel increases, making optimization more important. Thereby, the matching is generally subject to economies of scale that are not incorporated in the CA method. With economies of scale, we refer to the observation that if the number of parcels and crowd-shippers increase by the same proportion, the percentage of demand delivered generally increases. Quantification of this effect is outside the scope of this work. As supply increases further, a situation is reached where all demand to many region can be served easily. As the continuum approximation incorporates the upper bound imposed by available demand, we observe that performance increases again for the highest supply levels. 
\\\\
As the maximum detour increases, the difference between the approximation and the dynamic assignment decreases and therefore the quality of the approximation improves. The reason for this is that for higher values of $\tau$, as crowd-shippers can be feasibly matched to more parcels and parcels can be feasibly matched to more crowd-shippers, there exists more flexibility in the matching. That is, various matchings may lead to similar objective values. For the number of hubs, there is no clear effect on the quality of the approximation. We note that the percentual difference is substantially higher for the single-hub case, but this has to do with the relatively low objective values for this case. If we consider the absolute differences, this effect is far less substantial.

\subsection{Comparison of CA Approach to Simulation-Optimization Approach}\label{subsec:comparison_approach}
\begin{figure}[H]
    \centering
    \hspace*{-25mm}
    \includegraphics[width = 1.2\textwidth]{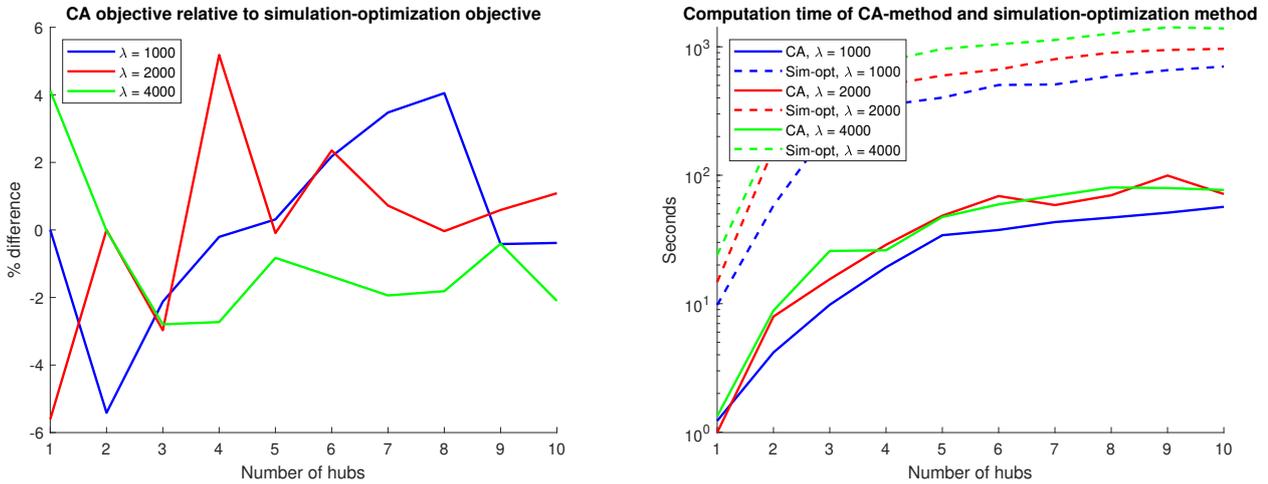}
    \caption{Comparison of CA and simulation-optimization method}
    \label{fig:computation_time}
\end{figure}

The previous section indicated the quality of the continuum approximation is high. In this section, we further evaluate the performance of our algorithm by comparing the CA-based solution algorithm proposed in this work to a well-known simulation-optimization approach \citep{carson1997simulation}. A simulation-optimization approach generally consists of two layers: an outer layer where an optimization problem is solved and an inner layer where a simulation is used to evaluate the performance of the optimized outer layer. This is a common approach for difficult (transportation) systems with many interactions, for example, \cite{romero2012simulation} and \cite{layeb2018simulation}. For the simulation-optimization approach, we use the same LNS algorithm but the continuum approximation is replaced by performing multiple simulations using the minimal detour assignment and taking the average number of delivered parcels. Here, for the sake of computation time, we use only 2 simulations. We note that further increasing the number of simulations has a proportional effect on the computation time, but can also lead to slightly better objectives. We compare the two methods in terms of objective value (i.e. number of parcels delivered by crowd-shippers) and computation time. The objective value is evaluated based on 10 simulation runs with a minimal-detour assignment strategy, that are different from the 2 simulations that were used at the optimization stage. 
\\\\
The results are displayed in Figure \ref{fig:computation_time}. The left-hand panel displays the objective of the CA method relative to the objective of the simulation-optimization method. A positive percentage indicates the CA method performs better whereas a negative percentage indicates that the simulation-optimization method performs better. Clearly, the objectives of the two methods are extremely close. More specifically, if we look at the exact opened hubs determined by each of the methods, these are also very similar. On the one hand, by increasing the number of simulations, the objective of the simulation-optimization method can be further improved. On the other hand, the CA-based dynamic assignment strategy can be leveraged to improve the dynamic assignment performance, as is discussed in Section \ref{subsec:comparison}. 
\\\\
The computation times of each of the methods are displayed in the right-hand panel of Figure \ref{fig:computation_time}, where we note the log-scale of the y-axis. 
We observe that the CA approach is significantly faster than the simulation-optimization approach. With the CA approach, we are able to determine the best set of hubs within approximately 1 minute, even for a high number of hubs. The computation time of the simulation-optimization method is between $10$ and $25$ times higher, depending on the supply level. We observe that the computation time of the CA algorithm is relatively more influenced by the number of hubs than the other two algorithms. However, even for a relatively high number of hubs, the computation time is still lower. In terms of the expected number of potential crowd-shippers, the computation time of the CA algorithm remains roughly the same. For the simulation-optimization method, the computation time increases proportionally to the number of supply units. 
\\\\
Clearly, the computation time of both methods also increases with the size of the network and the number of hubs that are needed. For larger networks, such as the complete city of Washington DC (rather than the smaller sample we have chosen here), the simulation-optimization approach will turn out to be too time-consuming to use. Especially in case the hubs are mobile and can be moved on a day-to-day basis, having a rapid way to determine their location is absolutely necessary. On top of this, the displayed computation times are for a single scenario where the demand level, supply level, number of hubs and prices are fixed. In reality, a large variety of scenarios may need to be tested to choose the most appropriate one, which drastically increases the total computation time and therefore the advantage of the CA approach.

\subsection{Results on the Network}\label{subsec:results_network}
\begin{figure}[H]
    \centering
      \hspace*{-7mm}
     \begin{subfigure}[b]{0.45\textwidth}  
            \centering
            \includegraphics[width=\textwidth]{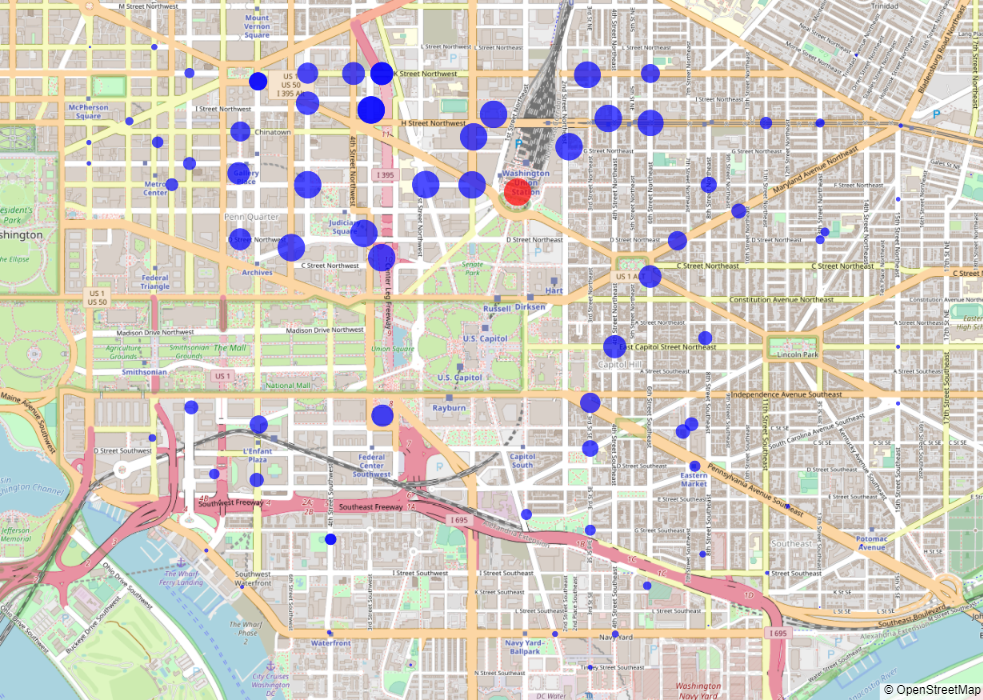}
            \caption{Single hub}
            \label{fig:network_1hub}
        \end{subfigure}
         \begin{subfigure}[b]{0.45\textwidth}  
            \centering
            \includegraphics[width=\textwidth]{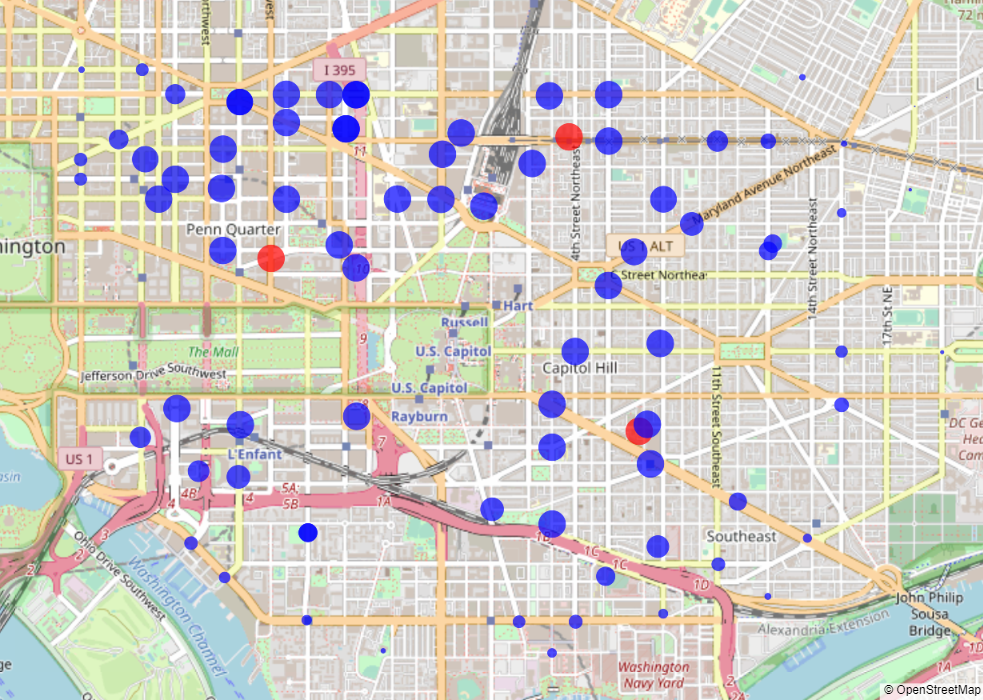}
            \caption{Three hubs}
            \label{fig:network_3hub}
        \end{subfigure}
    \caption{Bubble chart of demand served by crowd-shippers in the network}
    \label{fig:networks}
\end{figure}

In order to obtain managerial insights regarding the exact location of hubs, we evaluate the hub locations in the network. The bubble chart in Figure \ref{fig:networks} displays the considered network where a blue bubble represents a regular demand region and a red bubble represents a demand region that was chosen to be a hub location. The size of the bubble represents the fraction of demand that was served by crowd-shippers. That is, a full bubble implies that all demand in the region is expected to be delivered by crowd-shippers, according to the continuum approximation, whereas a smaller bubble implies that only a fraction of the demand is expected to be delivered by crowd-shippers. Figure \ref{fig:network_1hub} displays the single-hub approximation and Figure \ref{fig:network_3hub} displays the case where 3 hubs are used.

\begin{figure}[H]
    \centering
      \hspace*{-7mm}
     \begin{subfigure}[b]{0.45\textwidth}  
            \centering
            \includegraphics[width=\textwidth]{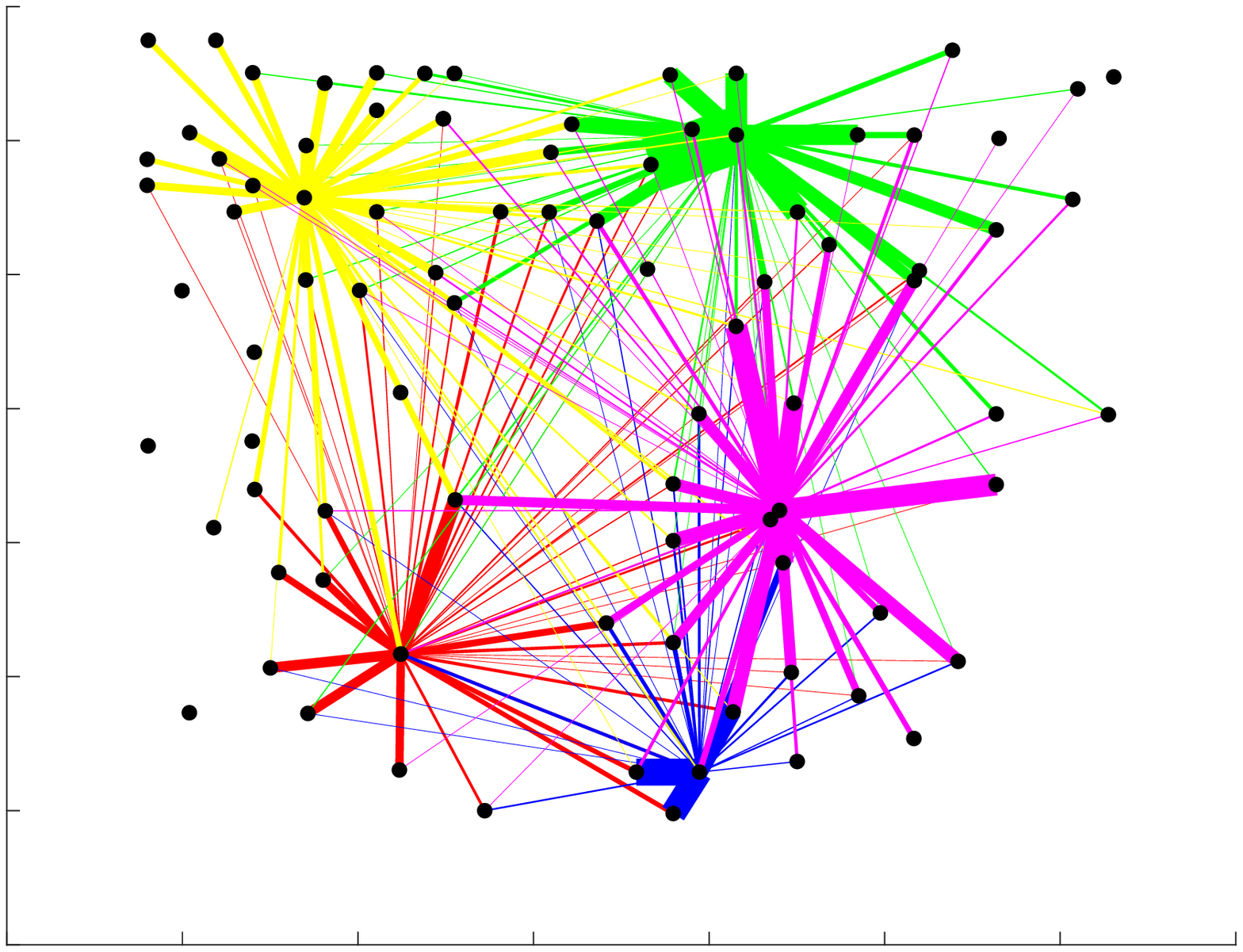}
            \caption{Spider chart}
            \label{fig:hubassignment}
        \end{subfigure}
         \begin{subfigure}[b]{0.45\textwidth}  
            \centering
            \includegraphics[width=\textwidth]{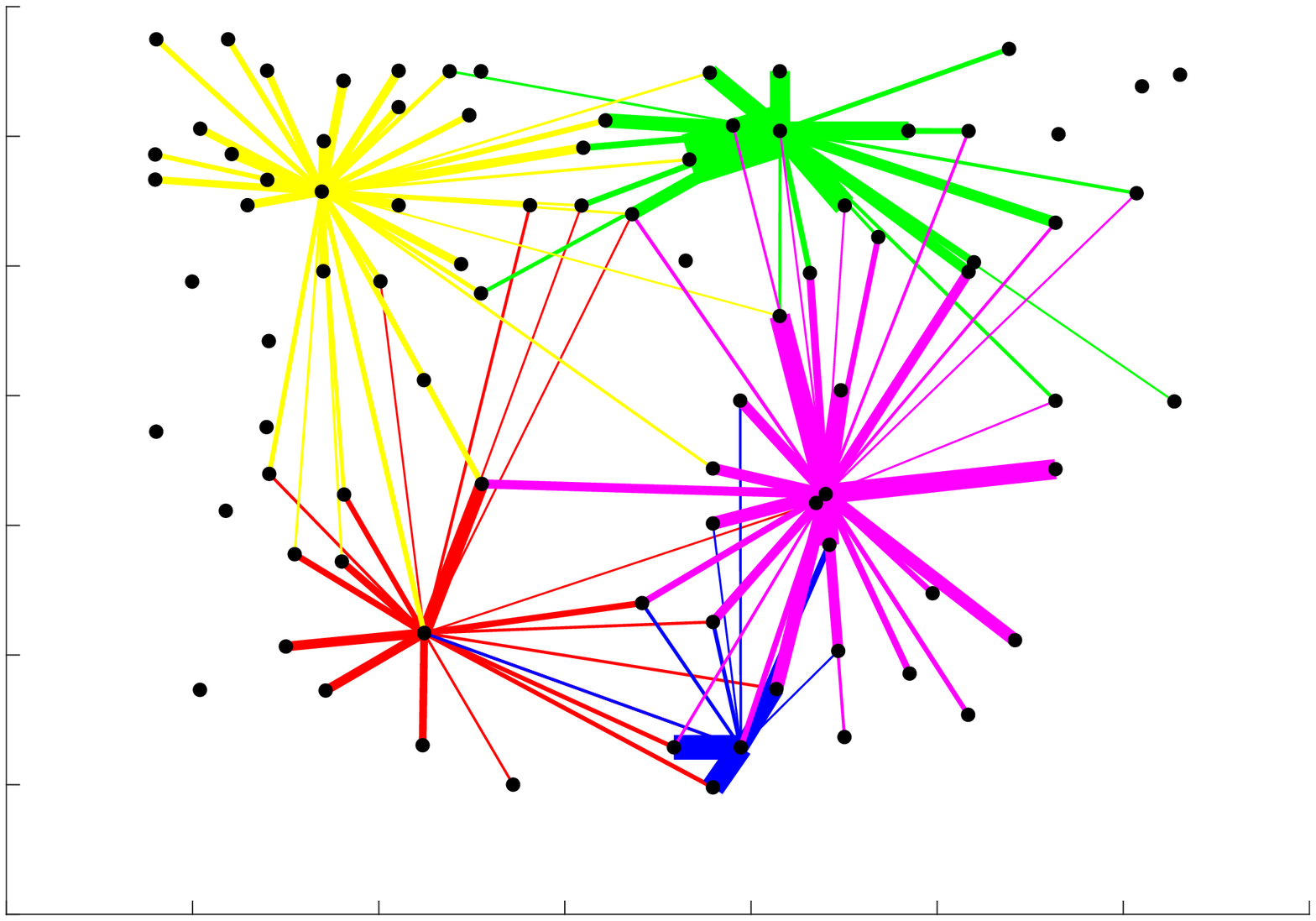}
            \caption{Truncated spider chart}
            \label{fig:hubassignment_trunc}
        \end{subfigure}
    \caption{Spider chart linking the origins and destinations of delivered parcels}
    \label{fig:spiders}
\end{figure}

\noindent When only a single hub is used, this hub is located at the train station. Horizontally, this location is central, but vertically this location is slightly shifted to the north. The main reason for this location is that the train station (Union Station) is one of the most popular origins in the region and therefore a substantial number of potential crowd-shippers can be attracted. In this case, we observe that the potential number of attracted crowd-shippers is more important than geographical centrality of the location. 
\\\\
When multiple hubs are built, we observe that the location near the train station is (approximately) preserved. On top of this, more locations are built that are spread throughout the city. Spreading the locations sufficiently throughout the city has two advantages. The first is that hubs serve different demand regions. By spreading the hubs, different demand regions are easier to reach by crowd-shippers and this therefore increases the number of served parcels. Secondly, different crowd-shippers (i.e. crowd-shippers with different origins and destinations) can be attracted. Therefore, spreading hubs sufficiently throughout the city allows to increase both the reachable demand and the attractable supply. 
\\\\
The strong inter-dependency between the hubs is shown in Figure \ref{fig:spiders}, which displays from which hubs the parcels delivered to each region originate. For this example, 5 hubs are used and the CA-based dynamic assignment strategy is implemented. Figure \ref{fig:hubassignment} shows the untruncated spider chart, where the line-width is determined by the number of parcels that have the same origin and destination. Figure \ref{fig:hubassignment_trunc} shows the truncated spider chart, where all lines with 3 or less parcels have been eliminated. We note that some regions are not served by crowd-shippers. This can have two reasons. The first, which holds for the regions at the edge of the network, is that they are actually not served by crowd-shippers, because there are few to no crowd-shippers that can feasibly reach that region. The second, which holds for some regions at the center of the network (corresponding to the Mall), is that demand is zero and therefore no service is required for those regions. 
\\\\
It is clear that the majority of the destination regions are served by multiple hubs. In Figure \ref{fig:hubassignment}, only 10\% of the destination regions are served by a single hub.After truncating in Figure \ref{fig:hubassignment_trunc}, still only 30\% of the regions are served by a single hub, the majority of which is caused by the yellow hub. The higher the number of hubs, the lower the number of regions served by a single hub and thus the higher the inter-dependency. Intuitively, this corresponds to the fact that when more hubs are constructed, the distance between hubs is lower and therefore their similarity (see Section \ref{subsec:stage_1}) increases. 

\subsection{Optimal Number of Hubs}\label{subsec:optimal_number_hubs}
\begin{figure}[H]
    \centering
    \includegraphics[width=1.0\textwidth]{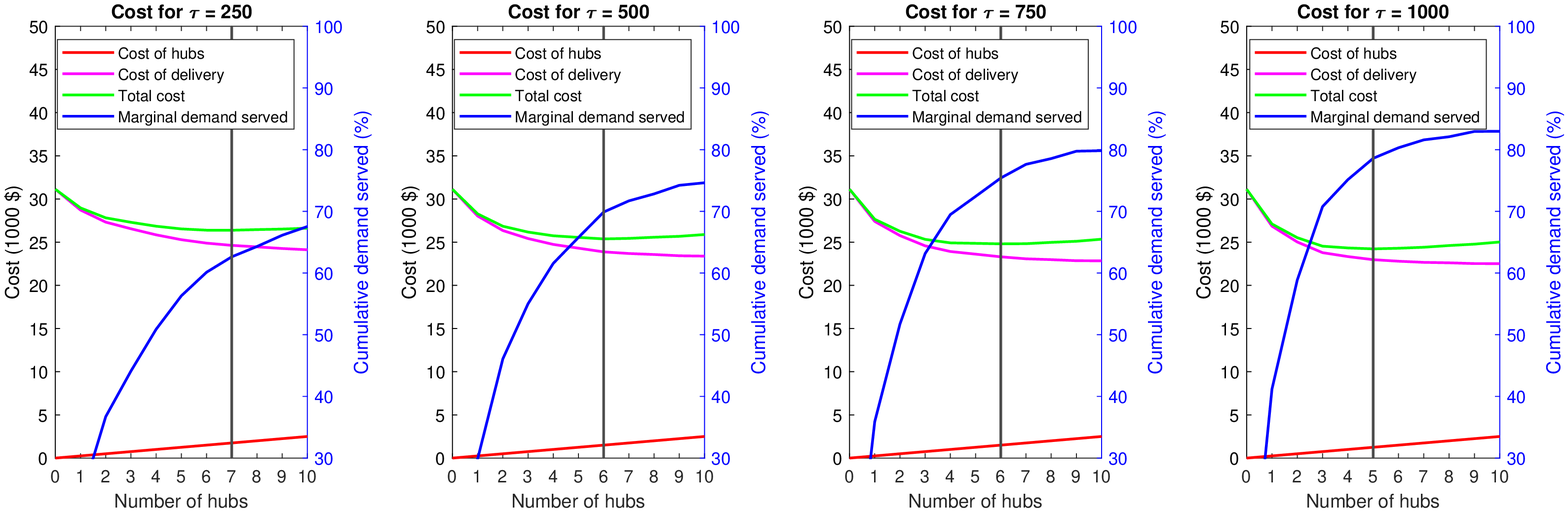}
    \caption{Decomposition of costs for varying number of hubs and maximum detour $\tau$}
    \label{fig:hubcost_tau}
\end{figure}

In this section, we evaluate the optimal number of hubs for various values of $\tau$. Figure \ref{fig:hubcost_tau} displays the costs for a varying number of hubs, as approximated by the CA approach. The costs are decomposed as delivery costs (the sum of regular delivery cost and crowd-shipping fees) and the fixed costs of operating hubs. Thereby, the blue curve displays the cumulative demand served as a percentage of the total demand of adding an additional hub.
\\\\
Increasing the value of $\tau$ has two opposing effects on the total costs. On the one hand, potential hub locations that previously were not able to serve sufficient demand to be profitable may now be able to serve more demand because the maximum detour is higher, thereby increasing the optimal number of hub locations. On the other hand, a hub may be able to reach more demand regions due to the increase of $\tau$, making other hubs obsolete, thereby reducing the optimal number of hubs. A similar effect is true for an increase or decrease of the total expected supply, which may either lead to an increase or decrease in the optimal number of hubs. 
\\\\
The marginal percentage of demand served due to the addition of one more hub is clearly diminishing (as shown in the right vertical axis of Figure \ref{fig:hubcost_tau}, where the curvature of the cumulative demand decreases). The first hubs are the most profitable and can serve a relatively large portion of the demand. As we showed in Section \ref{subsec:results_network}, these hubs are built at central locations where the number of potential crowd-shippers is high. Afterwards, additional hubs can be opened at less busy locations to further increase served demand, but the effect is substantially lower. 

\FloatBarrier
\FloatBarrier

\subsection{Comparison to non-predictive strategies}\label{subsec:parcel_hub_results}
As described in Section \ref{sec:methodology}, the proposed solution approaches to all three stages contain a predictive component in order to incorporate the influence of demand, supply and their interaction in the decision making. In this section, we compare CA-based strategies to strategies that do not use such a predictive component and only base the decisions based on geographical distances. The first stage hub locations are determined using a facility location problem minimizing the total distance of every region to the closest hub. The second stage parcel-hub assignment is made by assigning every parcel to the hub closest to its destination and the third stage matching is based on minimal-detour for crowd-shippers. Such a non-predictive approach ignores the connection between the three stages and solves every stage separately, basing its decisions solely on distance rather than demand and supply estimates and their interaction. 
\\\\
The results are displayed in Table \ref{tab:nonpredictive}. For every scenario, the number of served parcels is an average over 10 simulation runs. We observe that using CA-based strategies for all three stages significantly improves the performance of the system. By using CA estimates to determine the hub locations, parcel-hub assignment and the matching of parcels to crowd-shippers, between 7\% to 31\% more parcels can be served by crowd-shippers compared to when non-predictive strategies are used. Clearly, as $\tau$ increases and crowd-shippers become more flexible, the improvement CA makes diminishes. For the number of hubs, this is not straightforward.

\begin{table}[H]
\caption{Comparison of CA strategies to non-predictive strategies}
\label{tab:nonpredictive}
\begin{threeparttable}

\begin{tabular}{ccrrr}
\toprule 
$\tau$ & $|H|$ & Non-predictive & CA   & \% Improvement \\ \midrule
250 & 3 & 1281           & 1644 & 28.3       \\
250 & 5 & 1556           & 2039 & 31.1       \\
250 & 7 & 1866           & 2295 & 23.0       \\
500 & 3 & 1839           & 2087 & 13.5       \\
500 & 5 & 2097           & 2499 & 19.2       \\
500 & 7 & 2370           & 2692 & 13.6       \\
750 & 3 & 2270           & 2446 & 7.8        \\
750 & 5 & 2433           & 2681 & 10.2       \\
750 & 7 & 2626           & 2899 & 10.4      \\ \bottomrule
\end{tabular}
 \begin{tablenotes}[para]
\footnotesize{ The first two columns give the properties of the instance. The maximum detour $\tau$ is given in meters and the number of hubs is denoted by $|H|$. The third and fourth column denote the number of parcels served by crowd-shippers for non-predictive and CA strategies respectively and the last column denotes the improvement CA makes over the non-predictive strategies as a percentage. 
}
\end{tablenotes}

\end{threeparttable}

\end{table}

\subsection{Comparison of assignment policies for endogenous supply}\label{subsec:assignment_policies}

Using a simplified version of the acceptance model explored by \cite{wicaksono2021market}, we define the flow of potential crowd-shippers as a function of maximum detour and reward. The reward of $5\$$ for a maximum detour of 500 meters is used as a baseline. For every additional 500 meters, supply is expected to decrease by 10\%. On the other hand, every additional dollar increases the supply by 5\%. Specifically, we consider four options for $\tau = [500, 1000, 1500, 2000]$ and three options for the reward $p_{\tup {cs}} = [3, 5, 7]$. The value of $\lambda(\tau, p_{\tup {cs}})$ follows directly from these values using the described rule. Tables \ref{tab:minimal_detour}, \ref{tab:batch} and \ref{tab:ca} display various statistics for the minimal detour, batch and CA-based assignment strategy. The optimal hub locations are taken from the CA-based hub location algorithm, such that the optimal number of hubs is the same for all three strategies. The number of served parcels, total cost and average detour made by crowd-shippers are obtained as an average of 20 simulation runs. 
\\\\
First of all, we observe that in general, in agreement with the results discussed in Section \ref{subsec:optimal_number_hubs}, the optimal number of hubs decreases with the available supply $\lambda$ and decreases (partially through the supply) with the maximum detour $\tau$. For higher values of $\lambda$ and $\tau$, the same number of parcels can be served using fewer hubs. In addition to this, higher rewards $p_{\tup {cs}}$ reduce the benefit of a marginal hub and therefore also decrease the optimal number of hubs. 
\\\\
The minimal-detour assignment strategy clearly leads to a substantially lower average detour compared to the other two policies. The average detour of the CA-based assignment strategy is slightly higher than for the batch assignment, which is a consequence of nature of the strategy. As the CA-based assignment strategy favors parcels that are less likely to be delivered, these often lead to a higher detour. For all policies, we observe that the average detour is substantially lower than the maximum detour crowd-shippers are willing to make. This suggests that the cost can be further reduced by using a pricing policy that adapts the reward to the actual detour, rather than the maximum detour. However, this raises the concern of truthfulness in the reporting behaviour of crowd-shippers which has been a concern for many incentive-based policies in mobility systems \citep{asghari2017line, stokkink2021predictive}.
\\\\
The success of the CA-based strategy is clear from the high number of parcels served by crowd-shippers and the low corresponding costs. The CA-based strategy outperforms the other two strategies for almost all scenarios in terms of these two aspects. The best $(\tau, p_{\tup {cs}})$-combination in terms of total cost is marked in grey for every assignment strategy. For this, we see that the total cost for the CA-based assignment strategy are roughly 550\$ ($\approx 2.7\%$)  lower than the batch assignment and approximately 700\$ ($\approx 3.4\%$) compared to the minimal detour strategy.

\begin{table}[H]
\caption{Statistics for Minimal-Detour Assignment Strategy}
\label{tab:minimal_detour}
\small
\begin{threeparttable}
\begin{tabular}{ccccccc}
\toprule
$\tau$ & $p_{\tup {cs}}$ & $\lambda(\tau, p_{\tup {cs}})$ & $|H^*|$ & Served parcels & Total cost & Avg. detour \\ \midrule 
500  & 3      & 3784                  & 8     & 2487.7         & 21690.58   & 170.6          \\
500  & 5      & 4232                  & 6     & 2610.5         & 25858.63   & 172.1          \\
500  & 7      & 4612                  & 2     & 1815.1         & 30477.45   & 199.2          \\ \midrule
1000 & 3      & 3374                  & 6     & 2619.0         & 20599.73   & 333.4          \\
1000 & 5      & 3784                  & 5     & 2752.7         & 25253.25   & 399.9          \\
1000 & 7      & 4232                  & 2     & 2249.2         & 30260.40    & 426.5          \\ \midrule 
\rowcolor[HTML]{D9D9D9} 
1500 & 3      & 2955                  & 4     & 2530.9         & 20496.18   & 544.7          \\
1500 & 5      & 3374                  & 4     & 2759.0         & 24987.50    & 592.3          \\
1500 & 7      & 3784                  & 2     & 2458.1         & 30155.93   & 693.6          \\ \midrule 
2000 & 3      & 2532                  & 4     & 2345.6         & 21330.02   & 650.6          \\
2000 & 5      & 2955                  & 4     & 2622.2         & 25329.63   & 673.1          \\
2000 & 7      & 3374                  & 2     & 2566.0         & 30102.03   & 845.4          \\ \bottomrule
\end{tabular}
 \begin{tablenotes}[para]
\footnotesize{The maximum detour $\tau$ is given in meters and the crowd-shipper reward is given in dollars. The supply $\lambda$ is computed directly as a function of $\tau$ and $p_{\tup {cs}}$ and the optimal hub locations $H^*$ are determined using the CA approach. The number of served parcels and total cost (in dollars) are daily averages and the average detour is given in meters.}
\end{tablenotes}
\end{threeparttable}
\end{table}

\begin{table}[H]
\caption{Statistics for Batch Assignment Strategy}
\label{tab:batch}
\small
\begin{threeparttable}
\begin{tabular}{ccccccc}
\toprule
$\tau$ & $p_{\tup {cs}}$ & $\lambda(\tau, p_{\tup {cs}})$ & $|H^*|$ & Served parcels & Total cost & Avg. detour \\ \midrule 
500  & 3      & 3784                  & 8     & 2548.5         & 21416.98   & 281.7          \\
500  & 5      & 4232                  & 6     & 2678.5         & 25688.88   & 284.5          \\
500  & 7      & 4612                  & 2     & 1813.9         & 30478.05   & 294.4          \\ \midrule
\rowcolor[HTML]{D9D9D9} 
1000 & 3      & 3374                  & 6     & 2675.6         & 20344.80    & 666.2          \\
1000 & 5      & 3784                  & 5     & 2817           & 25092.50    & 691.2          \\
1000 & 7      & 4232                  & 2     & 2260           & 30255.00      & 692.1          \\ \midrule
1500 & 3      & 2955                  & 4     & 2546.8         & 20424.40    & 1040           \\
1500 & 5      & 3374                  & 4     & 2803           & 24877.63   & 1051.1         \\
1500 & 7      & 3784                  & 2     & 2528.1         & 30120.98   & 1120.1         \\ \midrule
2000 & 3      & 2532                  & 4     & 2399.4         & 21087.70    & 1359.9         \\
2000 & 5      & 2955                  & 4     & 2702           & 25130.13   & 1383.1         \\
2000 & 7      & 3374                  & 2     & 2616.2         & 30076.90    & 1463.7           \\ \bottomrule
\end{tabular}
 \begin{tablenotes}[para]
\footnotesize{The maximum detour $\tau$ is given in meters and the crowd-shipper reward is given in dollars. The supply $\lambda$ is computed directly as a function of $\tau$ and $p_{\tup {cs}}$ and the optimal hub locations $H^*$ are determined using the CA approach. The number of served parcels and total cost (in dollars) are daily averages and the average detour is given in meters.}
\end{tablenotes}
\end{threeparttable}
\end{table}

\begin{table}[H]
\caption{Statistics for CA-based Assignment Strategy}
\label{tab:ca}
\small
\begin{threeparttable}
\begin{tabular}{ccccccc}
\toprule
$\tau$ & $p_{\tup {cs}}$ & $\lambda(\tau, p_{\tup {cs}})$ & $|H^*|$ & Served parcels & Total cost & Avg. detour \\ \midrule 
500  & 3      & 3784                  & 8     & 2555.5         & 21385.48   & 288.3          \\
500  & 5      & 4232                  & 6     & 2649           & 25762.38   & 281            \\
500  & 7      & 4612                  & 2     & 1884.1         & 30442.95   & 303            \\ \midrule
1000 & 3      & 3374                  & 6     & 2752.1         & 20000.55   & 691.7          \\
1000 & 5      & 3784                  & 5     & 2870.1         & 24959.75   & 701.9          \\
1000 & 7      & 4232                  & 2     & 2363.1         & 30203.48   & 716            \\ \midrule
\rowcolor[HTML]{D9D9D9} 
1500 & 3      & 2955                  & 4     & 2687.1         & 19793.05   & 1133.7         \\
1500 & 5      & 3374                  & 4     & 2924.8         & 24573      & 1120.4         \\
1500 & 7      & 3784                  & 2     & 2674.8         & 30047.63   & 1184.1         \\ \midrule
2000 & 3      & 2532                  & 4     & 2492.9         & 20667.18   & 1572.3         \\
2000 & 5      & 2955                  & 4     & 2849.1         & 24762.25   & 1610.8         \\
2000 & 7      & 3374                  & 2     & 2801           & 29984.53   & 1635.7           \\ \bottomrule
\end{tabular}
 \begin{tablenotes}[para]
\footnotesize{The maximum detour $\tau$ is given in meters and the crowd-shipper reward is given in dollars. The supply $\lambda$ is computed directly as a function of $\tau$ and $p_{\tup {cs}}$ and the optimal hub locations $H^*$ are determined using the CA approach. The number of served parcels and total cost (in dollars) are daily averages and the average detour is given in meters.}
\end{tablenotes}
\end{threeparttable}
\end{table}

\noindent Using an adaptive pricing mechanism can form a large improvement over the fixed rewards considered here. In such a pricing scheme, higher rewards are given for high-demand and low-supply routes, whereas lower rewards are given for routes that are relatively popular among potential crowd-shippers. Clearly, this is outside the scope of this work, but it marks an important direction of future research.

\section{Conclusion}\label{sec:conclusion}
Crowd-shipping is a promising alternative to traditional last-mile delivery methods that can help to reduce congestion, pollution and improve the overall performance of the delivery system. One of the main drivers is the availability of sufficient suppliers. In this paper, we proposed a hub-based crowd-shipping system where crowd-shippers pickup parcels from hub locations and deliver them to the final destination. By constructing hubs at strategic locations, more crowd-shippers can be attracted to serve more demand. 
\\\\
We developed a heuristic approach to efficiently determine the optimal hub locations. The operational problem of assigning parcels to potential crowd-shippers is approximated using a Continuum Approximation (CA) approach. Using this, the best set of hubs is explored using a large neighborhood search heuristic. Using quality and similarity metrics, the search space is explored efficiently and a solution is found within a reasonable amount of time.
\\\\
The approximation shows to provide an accurate estimation of the number of served parcels mostly being within 20\% of the optimal static assignment problem and within 12\% of the more realistic dynamic assignment problem. Performance is clearly higher when crowd-shippers are willing to make larger detours. Although performance decreases when potential supply increases, using other methods becomes nearly impossible for high levels of supply due to the substantial increase in computation time. Other than simulation or optimization approaches, the computation time of the continuum approximation approach is independent of the expected number of potential crowd-shippers. 
\\\\
We designed a dynamic assignment strategy based on the obtained CA estimates. Using a dynamic discrete event simulator, the CA-based assignment strategy is compared to a minimal-detour assignment strategy as well as a batch assignment strategy. An analysis over a broad set of parameters shows that the CA-based assignment strategy generally outperforms the other strategies in terms of the number of parcels that can be served by crowd-shippers, as well as the total system cost. On top of this, our results show that it is not necessarily best to assign a parcel to the closest hub in terms of distance and significant improvements may be achieved by using more advanced parcel-hub assignment strategies. We show that our CA-based strategies for the three-stage problem outperforms a strategy where the decisions in all three stages are made based on non-predictive methods. Our approach performs up to 31\% better than non-predictive methods on the evaluated scenarios. 
\\\\
A natural extension of this work is to incorporate vehicle routing costs directly in the hub location problem. The framework designed in this paper allows to efficiently incorporate routing decisions as long as the problem is solved or approximated fast. In addition to this, considering crowd-shippers that carry multiple parcels is an interesting direction of future research. Furthermore, a generalization of the problem to include (soft) time-windows would enable customers to indicate their preferred delivery time. A dynamic pricing policy for heterogeneous population in terms of value of time, maximum detour and other parameters would also shed some light in the potential of real implementations of these strategies.


\color{black}
\clearpage

\bibliography{references.bib}
\bibliographystyle{apalike}

\end{document}